\documentclass[12pt,a4paper]{amsart}
\usepackage{geometry}
\usepackage{tabls}
\usepackage{url}
\usepackage[utf8]{inputenc}
\usepackage{amsfonts}
\usepackage{amssymb}
\usepackage{mathrsfs}
\usepackage{amsmath}
\usepackage{amsthm}
\usepackage{amscd}
\usepackage{fancyhdr}
\usepackage{enumerate}
\usepackage{paralist}
\usepackage{xypic}
\usepackage{graphicx}
\usepackage{cite}
\usepackage{lipsum}
\usepackage{footmisc}
\usepackage[all,cmtip]{xy}

\numberwithin{equation}{section}

\theoremstyle{plain}
\newtheorem{Cor}{Corollary}[section]
\newtheorem{Def}[equation]{Definition}
\newtheorem{Thm}[equation]{Theorem}
\newtheorem{lem}[equation]{Lemma}
\newtheorem{prop}[equation]{Proposition}
\newtheorem{rem}[equation]{Remark}

\newtheorem{conj}[equation]{Conjecture}

\begin{document}

\title{The depth structure of motivic multiple zeta values}

\author{Jiangtao Li}

\address{School of Mathematic Sciences\\
        Peking University\\
         Beijing, China}

\email{ljt-math@pku.edu.cn}

\begin{abstract}

   In this paper, we construct some maps related to the motivic Galois action on depth-graded motivic multiple zeta values. From these maps we give some short exact sequences about depth-graded motivic multiple zeta values in depth two and three. In higher depth we conjecture that there are exact sequences of the same type. We will show from three conjectures about depth-graded motivic Lie algebra we can nearly deduce the exact sequences conjectures in higher depth. At last we give a new proof of  the  result that the modulo $\zeta^{\mathfrak{m}}(2)$ version motivic double zeta values are generated by the totally odd part. We
   reduce the well-known conjecture that the modulo $\zeta^{\mathfrak{m}}(2)$ version motivic triple zeta values are generated by the totally odd part
   to an isomorphism conjecture in linear algebra.
\end{abstract}
\maketitle
\let\thefootnote\relax\footnotetext{
2010 $\mathnormal{Mathematics} \;\mathnormal{Subject}\;\mathnormal{Classification}$. Primary 11F32, Secondary 11F67.\\
$\mathnormal{Keywords:}$  Multiple zeta values, Period polynomial, Lie algebra, Motives. }
\section{Introduction}\label{int}

The multiple zeta values are defined by convergent series
\[\zeta(n_1,...,n_r)=\sum_{0<k_1<\cdots<k_r}\frac{1}{k_1^{n_1}\cdots k_r^{n_r}}, \, (n_1,...,n_{r-1}>0,n_r>1).\]
Multiple zeta values  play important roles in many fields. For $\zeta(n_1,..., n_r)$,  $N=n_1+\cdots+n_r$ and $r$ are called its weight and depth respectively. Let $\mathcal{Z}_0=\mathbb{Q}$ and  $\mathcal{Z}_N$ the space of $\mathbb{Q}$-linear combinations of multiple zeta values of weight $N$. Then it's clear that
\[
\mathcal{Z}=\bigoplus_{N=0}^{+\infty}\mathcal{Z}_N
\]
is a graded commutative algebra over $\mathbb{Q}$.

In \cite{brown}, Brown defined the motivic multiple zeta
 values by using an idea of Goncharov \cite{goncha}. Motivic multiple zeta values  satisfy the double shuffle relations \cite{ihara},\cite{racinet}, \cite{souderes},  and there's a motivic Galois action on them.
Their elements are $\mathbb{Q}$-linear combinations of motivic multiple zeta values $\zeta^{\mathfrak{m}}(n_1,...,n_r)$.
What's more, there is a surjective graded algebra homomorphism:
\[
\eta:\mathcal{H}\rightarrow\mathcal{Z}
\]
\[
\zeta^{\mathfrak{m}}(n_1,...,n_r)\rightarrow\zeta(n_1,...,n_r)
\]

 In \cite{brown}, Brown proved a conjecture of Hoffman (Conjecture C in \cite{hoff}) in multiple zeta values by using the motivic Galois action on motivic multiple zeta values. So the study of motivic multiple zeta values is helpful to understand the structure of multiple zeta values.

There is a depth filtration on $\mathcal{H}$. For motivic multiple zeta value $\zeta^{\mathfrak{m}}(n_1,...,n_r)$, we call $r$ its depth. Denote by $\mathfrak{D}_r\mathcal{H}$ the space of $\mathbb{Q}$-linear
combinations of motivic multiple zeta values of depth $\leq r$. Define \[gr^{\mathfrak{D}}_r\mathcal{H}=\mathfrak{D}_r\mathcal{H}/\mathfrak{D}_{r-1}\mathcal{H}.\]

In \cite{kreimer}, Broadhurst and Kreimer proposed a conjecture about the depth structure of multiple zeta values. In the motivic setting, Brown proposed the following motivic Broadhurst-Kreimer conjecture in \cite{depth}.
\begin{conj}(Broadhurst-Kreimer-Brown)\label{bkb}
Denote by $gr_r^{\mathfrak{D}}\mathcal{H}_{N}$ the weight $N$ part of $gr_r^{\mathfrak{D}}\mathcal{H}$, then
\[
1+\sum_{N,r>0}\mathrm{dim}_{\mathbb{Q}}(gr_r^{\mathfrak{D}}\mathcal{H}_N)x^Ny^r=\frac{1+\mathbb{E}(x)y}{1-\mathbb{O}(x)y+\mathbb{S}(x)y^2-\mathbb{S}(x)y^4},
\]
where
\[
\mathbb{E}(x)=\frac{x^2}{1-x^2},\;\mathbb{O}(x)=\frac{x^3}{1-x^2},\;\mathbb{S}(x)=\frac{x^{12}}{(1-x^4)(1-x^6)}.
\]
\end{conj}

The case $r=2$ is a consequence of  the work of Gangl, Kaneko, Zagier \cite{kaneko}, a special case of the parity results about double shuffle relations which was firstly discovered by Euler and some basic facts about motivic multiple zeta values. The case $r=3$ is known due to the work of Goncharov \cite{goncharov}, \cite{gon}. The higher depth cases remain unknown.

If $n_1,...,n_r\geq 3$ and all of them are odd numbers, then we call the motivic multiple zeta value $\zeta^{\mathfrak{m}}(n_1,...,n_r)$ totally odd. Denote by $gr^{\mathfrak{D}}_r\mathcal{A}^{odd}$ the image of $\mathbb{Q}$-linear combinations of totally odd motivic multiple zeta values of depth $\leq r$ in $gr^{\mathfrak{D}}_r\mathcal{H}$ modulo $(\zeta^{\mathfrak{m}}(2))\cap gr^{\mathfrak{D}}_r\mathcal{H}$. In \cite{depth}, Brown proposed the following uneven part of motivic Broadhurst-Kreimer conjecture.
\begin{conj}\label{un}(Brown)
Denote by $gr_r^{\mathfrak{D}}\mathcal{A}^{odd}_N$ the weight $N$ part of $gr_r^{\mathfrak{D}}\mathcal{A}_N$, then
\[
1+\sum_{N,r>0}\mathrm{dim}_{\mathbb{Q}}(gr_r^{\mathfrak{D}}\mathcal{A}^{odd}_N)x^N y^r=\frac{1}{1-\mathbb{O}(x)y+\mathbb{S}(x)y^2}.
\]
\end{conj}

In Conjecture \ref{bkb} and Conjecture \ref{un}, the series $\mathbb{E}(x)$ and $\mathbb{O}(x)$ have geometric explanations that they are the generating series of the dimension space of even and odd single zeta values respectively. While $\mathbb{S}(x)$ is the generating series of the dimension space of the cusp forms of $SL_2(\mathbb{Z})$.

Since each term in Conjecture \ref{bkb} and Conjecture \ref{un} has a geometric explanation, we believe that there are some inner structures behind Conjecture \ref{bkb} and Conjecture \ref{un}. The main purpose of this paper is trying to understand these structures behind Conjecture \ref{bkb} and Conjecture \ref{un}.

   In fact, we have
\begin{Thm}\label{odd}
For $r=2$, there is an exact sequence
\[
0\rightarrow gr^{\mathfrak{D}}_2\mathcal{A}^{odd}\xrightarrow{\widetilde{\partial}} gr^{\mathfrak{D}}_1\mathcal{A}^{odd}\otimes gr^{\mathfrak{D}}_1\mathcal{A}^{odd}\xrightarrow{v} \mathbb{P}^{{}^{\vee}}\rightarrow 0.
\]
\end{Thm}

In the above short exact sequence, $\mathbb{P}$ denotes the restricted even period polynomial space which will be explained in detail later,  $\mathbb{P}^{{}^{\vee}}$ denotes its weight dual vector space (taking dual in each weight). The definitions of the maps involved will be given in the next section.

The  injectivity of $\widetilde{\partial}$ follows from taking dual. The relation between the free lie algebra  generated by two words with Ihara bracket and even restricted period polynomials was first discovered by Ihara and Takao  \cite{ihara1}. The exactness of the middle part follows from the main result of Schneps  \cite{sch} and the fact $gr_2^{\mathfrak{D}}\mathcal{A}^{odd}=gr_2^{\mathfrak{D}}\mathcal{A}$. So the exact sequence we have established provides an explanation of the relation between the work of Gangl, Kaneko, Zagier \cite{kaneko} and Schneps \cite{sch}.

In \cite{ll}, Fei Liu and the author established a short exact sequence for the subspace generated by the images of $\zeta^{\mathfrak{m}}(r,N-r)$, where $3\leq r\leq N-2$, odd in $gr_2^{\mathfrak{D}}\mathcal{H}_N$ for $N\geq 5$, odd. Theorem $1.1$ in \cite{ll} is  similar to Theorem \ref{odd} but the term which involves period polynomial becomes more complicated.

Inspired by Theorem \ref{odd}, for $r\geq 3$, we construct the maps $\widetilde{\partial}:gr^{\mathfrak{D}}_r\mathcal{A}^{odd}\rightarrow gr^{\mathfrak{D}}_{1}\mathcal{A}^{odd}\otimes gr^{\mathfrak{D}}_{r-1}\mathcal{A}^{odd}$ and  $D: gr^{\mathfrak{D}}_{1}\mathcal{A}^{odd}\otimes gr^{\mathfrak{D}}_{r-1}\mathcal{A}^{odd}\rightarrow \mathbb{P}^{{}^{\vee}}\otimes gr^{\mathfrak{D}}_{r-2}\mathcal{A}^{odd}$. We propose the following conjecture

\begin{conj}\label{ttt}
For $r\geq3$, there is an exact sequence

\[
 0\rightarrow gr^{\mathfrak{D}}_r\mathcal{A}^{odd}\xrightarrow{\widetilde{\partial}} gr^{\mathfrak{D}}_{1}\mathcal{A}^{odd}\otimes gr^{\mathfrak{D}}_{r-1}\mathcal{A}^{odd}\xrightarrow{D} \mathbb{P}^{{}^{\vee}}\otimes gr^{\mathfrak{D}}_{r-2}\mathcal{A}^{odd}\rightarrow 0.
\]

\end{conj}

We will show that why this conjecture should be true by analyzing the structure of the depth-graded motivic Lie algebra. We will prove that

\begin{Thm}\label{rrr}

(i)For $r=3$, the map $\widetilde{\partial}$ is injective and the map $D$ is surjective.
(ii)For $r\geq3$, assuming three conjectures about depth-graded motivic lie algebra, there will be an exact sequence
\[
 0\rightarrow gr^{\mathfrak{D}}_r\mathcal{A}^{odd}\xrightarrow{\widetilde{\partial}} gr^{\mathfrak{D}}_{1}\mathcal{A}^{odd}\otimes gr^{\mathfrak{D}}_{r-1}\mathcal{A}^{odd}\xrightarrow{D} \mathbb{P}^{{}^{\vee}}\otimes gr^{\mathfrak{D}}_{r-2}\mathcal{A}^{odd}\rightarrow 0.
\]

\end{Thm}

 The three conjectures (non-degenerated conjecture, vanishing conjecture and surjective conjecture) involved will be given in Section \ref{nonde}.

What's more we will show the relationship between the results of Tasaka and Theorem \ref{rrr}. Thus we provide a geometric explanation of the linear map constructed by Tasaka in \cite{tasaka}.

Denote by $\mathcal{A}=\mathcal{H}/(\zeta^{m}(2))$ and \[gr_r^{\mathfrak{D}}\mathcal{A}=gr_r^{\mathfrak{D}}\mathcal{H}/(\zeta^{\mathfrak{m}}(2))\cap gr_r^{\mathfrak{D}}\mathcal{H}\]
we can construct the map $\widetilde{\partial}_g:gr^{\mathfrak{D}}_r\mathcal{A}\rightarrow gr^{\mathfrak{D}}_{1}\mathcal{A}^{odd}\otimes gr^{\mathfrak{D}}_{r-1}\mathcal{A}$ and the map $D_g: gr^{\mathfrak{D}}_{1}\mathcal{A}^{odd}\otimes gr^{\mathfrak{D}}_{r-1}\mathcal{A}\rightarrow \mathbb{P}^{{}^{\vee}}\otimes gr^{\mathfrak{D}}_{r-2}\mathcal{A}$ the same way as $\widetilde{\partial}$ and $D$ respectively.

Then we have the following theorem
\begin{Thm}\label{gthree}
(i)For $r=2$, there is an exact sequence

\[0\rightarrow gr^{\mathfrak{D}}_2\mathcal{A}\xrightarrow{\widetilde{\partial}_g} gr^{\mathfrak{D}}_1\mathcal{A}^{odd}\otimes gr^{\mathfrak{D}}_1\mathcal{A}\xrightarrow{v_g} \mathbb{P}^{{}^{\vee}}\rightarrow 0.\]

(ii)For $r=3$, there is an exact sequence
\[
 0\rightarrow gr^{\mathfrak{D}}_3\mathcal{A}\xrightarrow{\widetilde{\partial}_g} gr^{\mathfrak{D}}_{1}\mathcal{A}^{odd}\otimes gr^{\mathfrak{D}}_{2}\mathcal{A}\xrightarrow{D_g} \mathbb{P}^{{}^{\vee}}\otimes gr^{\mathfrak{D}}_{1}\mathcal{A}\rightarrow 0.
\]
\end{Thm}
It's easy to show  that $gr_r^{\mathfrak{D}}\mathcal{A}^{odd}=gr_r^{\mathfrak{D}}\mathcal{A}$ for $r=1,2$ from the results of \cite{brown},\cite{kaneko},\cite{souderes}. The motivic Broadhurst-Kreimer conjecture (Conjecture \ref{bkb}) and its uneven part conjecture (Conjecture \ref{un}) suggest that $gr_r^{\mathfrak{D}}\mathcal{A}^{odd}=gr_r^{\mathfrak{D}}\mathcal{A}$ should also hold for $r=3$.

The proof of Theorem \ref{gthree} $(ii)$ is  mainly from the result of Concharov \cite{gon}.

Furthermore we propose the following conjecture
which shows the existence of a long exact sequence in general case:
\begin{conj}\label{lon}
For $r\geq4$, there is an exact sequence
\[
0\rightarrow\mathbb{P}^{{}^{\vee}}\otimes gr_{r-4}^{\mathfrak{D}}\mathcal{A}\rightarrow gr_r^{\mathfrak{D}}\mathcal{A}\rightarrow gr_1^{\mathfrak{D}}\mathcal{A}^{odd}\otimes gr_{r-1}^{\mathfrak{D}}\mathcal{A}\rightarrow \mathbb{P}^{{}^{\vee}}\otimes gr_{r-2}^{\mathfrak{D}}\mathcal{A}\rightarrow0.
\]

\end{conj}

The first map is not clear enough at present. It should have connection with the exceptional elements which are constructed by Brown in\cite{depth} if the exceptional elements are motivic.

The following theorem shows the importance of the non-degenerated conjecture to tackle Conjecture \ref{lon}.
\begin{Thm}\label{rere}
Assuming the non-degenerated conjecture for all $r\geq3$, then there is an exact sequence for each $r\geq3$.
\[
gr_r^{\mathfrak{D}}\mathcal{A}\xrightarrow{\widetilde{\partial}_g} gr_1^{\mathfrak{D}}\mathcal{A}^{odd}\otimes gr_{r-1}^{\mathfrak{D}}\mathcal{A}\xrightarrow{D_g} \mathbb{P}^{{}^{\vee}}\otimes gr_{r-2}^{\mathfrak{D}}\mathcal{A}\rightarrow0
\]
\end{Thm}

From Conjecture \ref{lon}  we can deduce the motivic Broadhurst-Kreimer conjecture easily. Furthermore, if Conjecture \ref{lon} is true, from the analysis of Enriquez and Lochak \cite{enriquez}, we will know the structure of the depth-graded motivic lie algebra and its  Lie algebra homology in each weight and depth.

The category equivalence between finite dimensional nilpotent lie algebra and unipotent algebraic group in characteristic 0 plays a key role in the proofs of our main results. So we give a short introduction in section \ref{key}.

We point out that if we can find an injective map $j:\mathbb{P}^{{}^{\vee}}\otimes gr_{r-4}^{\mathfrak{D}}\mathcal{A}\rightarrow gr_r^{\mathfrak{D}}\mathcal{A}$ for each $r\geq4$ such that
$\mathrm{Im}\,j\subseteq \mathrm{Ker}\,\widetilde{\partial}$, then Conjecture \ref{lon} will follows from the non-degenerated conjecture since the dimension of $\mathcal{A}$ in each weight is known by the work of Brown \cite{brown}.

The last section of this paper is devoted to the study of totally odd motivic multiple zeta values in depth $2$ and $3$. Using the work of Baumard and Schneps \cite{schneps} we prove that $gr_2^{\mathfrak{D}}\mathcal{A}^{odd}=gr_2^{\mathfrak{D}}\mathcal{A}$. This statement can also be proved by using the fact that motivic double zeta values satisfy double shuffle relations and main theorems in \cite{kaneko}.  In depth three, we show that a surjective conjecture proposed by Tasaka in linear algebra is suffice to deduce $gr_3^{\mathfrak{D}}\mathcal{A}^{odd}=gr_3^{\mathfrak{D}}\mathcal{A}$.

In this paper all the dual vector space means compact dual (linear functional whose support is finite dimension). Since motivic multiple zeta values satisfy double shuffle relations \cite{ihara},\cite{racinet}, \cite{souderes}.
Any statement holds for the classical multiple zeta values which is deduced only by double shuffle relations  also holds for the motivic multiple zeta values. 

\section{Motivic Galois action}\label{galois}

   In this section we will define all the maps involved in the short exact sequence in section \ref{int}.
These maps essentially come from the Galois action of motivic fundamental group of mixed Tate motives. The main references of this section are \cite{depth}, \cite{BF}, \cite{deligne}.

Denote by $\mathcal{MT}(\mathbb{Z})$ the category of mixed Tate motives over $\mathbb{Z}$. Let $\pi_1(\mathcal{MT}(\mathbb{Z}))$ its de-Rham fundamental group. From \cite{deligne}, we have 
\[
  \pi_1(\mathcal{MT}(\mathbb{Z}))=\mathbb{G}_m \ltimes \mathnormal{U}^{dR}.
  \]
Where the subgroup $\mathnormal{U}^{dR}$ is pro-unipotent and its graded Lie algebra (respective to the action of $\mathbb{G}_m$) is isomorphic to the free graded Lie algebra with one generator $\sigma_{2n+1}$ in every negative degree $-(2n+1),\; n\geq 1$.

Let ${}_0\Pi_1$ be the De-Rham realization of the motivic fundamental groupoid from
$\overrightarrow{1}_0$ to $\overrightarrow{-1}_1$ (the tangent vector $\overrightarrow{1}$ at $0$ and the tangent vector $\overrightarrow{-1}$ at $1$). Its function ring over $\mathbb{Q}$ is isomorphic to
\[\mathcal{O}({}_0\Pi_1)\simeq\mathbb{Q}\langle e^0,e^1\rangle,\]
where $\mathbb{Q}\langle e^0,e^1\rangle$ is equipped with the shuffle product. $\mathcal{O}({}_0\Pi_1)$ has a coaction from the function ring $\mathcal{O}(\mathrm{U}^{dR})$.

Denote by $dch\in {}_0\Pi_1(\mathbb{R})$ the de-Rham image of the straight line from $0$ to $1$. 
It determines a homomorphism
\[dch:\mathcal{O}({}_0\Pi_1)\longrightarrow\mathbb{R}.\]
Let $\omega_{e^0}(t)=\frac{dt}{t},\omega_{e^1}(t)=\frac{dt}{1-t}$. For any word $u_1\cdots u_k$ in $e^0,e^1$, if $\epsilon,\eta \to 0$, it's easy to check that 
$$\int\limits_{\varepsilon<t_{1}<\cdots <t_{k}<1-\eta}\omega_{u_{1}}(t_{1})\cdots \omega_{u_{k}}(t_{k})=P(log(\varepsilon), log(\eta))+O\left(\mathrm{sup}(\varepsilon|log(\epsilon)|^{A}+\eta|log(\eta)|^{B})\right),$$
where $P$ is a polynomial. Then we have $dch(u_1\cdots u_k)=P(0,0)$.

Let $I\subseteq\mathcal{O}({}_0\Pi_1)$ be the kernel of $dch$. Denote by $J^{\mathcal{MT}}\subseteq I$ be the largest graded ideal contained in $I$ which is stable under the coaction of $\mathcal{O}(\mathrm{U}^{dR})$.

Define $\mathcal{H}=\mathcal{O}({}_0\Pi_1)/J^{\mathcal{MT}}$ and the image of the word
\[
e^1\underbrace{e^0...e}_{n_1-1}{}^0e^1...e^1\underbrace{e^0...e}_{n_r-1}{}^0
\]
as $\zeta^{\mathfrak{m}}(n_1,...,n_r)$.

From the work of Brown \cite{brown}, we know that the elements of $\mathcal{H}$ are $\mathbb{Q}$-linear combinations of $\zeta^{\mathfrak{m}}(n_1,...,n_r)$. The motivic multiple zeta values satisfy the usual double shuffle relations \cite{ihara}, \cite{racinet}, \cite{souderes}.

Denote by ${}_x\Pi_y$ the de-Rham realization of motivic torsor of paths on $\mathbb{P}^1\backslash\{0,1,\infty\}$ from $x$ to $y$ where $x,y\in \{\overrightarrow{1}_0,\overrightarrow{-1}_1\}$. For convenience we write $\overrightarrow{1}_0$, $\overrightarrow{-1}_1$ as $0$, $1$ respectively. Denote by $G$
 the group of automorphisms of the groupoid ${}_x\Pi_y$ for $x,y\in{0,1}$ which respect to the following structures:

(1) (Groupoid structure) The multiplication maps
\[
   {}_x\Pi_y \times {}_y\Pi_z\rightarrow {}_x\Pi_z
   \]
for all $x,y,z\in\{0,1\}$.

(2) (Inertia) The automorphism fixes the elements
\[\mathrm{exp}(e_0)\in {}_0\Pi_0(\mathbb{Q}),\, \mathrm{exp}(e_1)\in {}_1\Pi_1(\mathbb{Q}),\]
where $e_0, e_1$ respectively denotes the differential $\frac{dz}{z},\frac{dz}{1-z}$.

Then by Proposition 5.11 in \cite{deligne}, we know that ${}_x\Pi_y$ is an $G$-torsor.
Since the groupoid and inertia structures are preserved by $\mathrm{U}^{dR}$, we have a morphism
\[
\varphi:\mathrm{U}^{dR}\rightarrow G\simeq{}_0\Pi_1.
\]

By the main results of \cite{brown}, $\varphi$ is injective. Denote by $\mathfrak{g}$ the corresponding Lie algebra of $\mathrm{U}^{dR}$, we have an injective map on Lie algebra
\[
i:\mathfrak{g}\rightarrow \mathrm{Lie}\, G \simeq(\mathbb{L} (e_0,e_1),\{ \; , \}).
\]
Where $\mathbb{L}(e_0,e_1)$ is isomorphic to the free Lie algebra generated by $e_0,e_1$ as a vector space, but its Lie algebra bracket is given by the Ihara bracket\cite{deligne},\cite{sch}.

Denote by $\mathfrak{h}=\left(\mathbb{L}(e_0,e_1),\{,\}\right)$ for convenience and we define a decreasing depth filtration on $\mathfrak{h}$ by
\[
\mathfrak{D}^r\mathfrak{h}=\{\xi\in \mathfrak{h}\mid \mathrm{deg}_{e_1}\;\xi\geq r\}.
\]
From Th$\acute{e}$or$\grave{e}$me 6.8(i) in\cite{deligne}, we know that
$i(\sigma_{2n+1})=(\mathrm{ad}\,e_0)^{2n}(e_1)+$ terms of degree $\geq 2 $ in $e_1$. So $\mathfrak{g}$ has an induced depth filtration. Since the Ihara bracket $\{\,,\,\}$ is compatible with the depth filtration,
\[\mathfrak{d}\mathfrak{g}=gr^{\bullet}\mathfrak{g}=\oplus_{r\geq1}\mathfrak{D}^r\mathfrak{g}/\mathfrak{D}^{r+1}\mathfrak{g}\]
is also a Lie algebra with the induced Ihara bracket.

Denote by $\mathfrak{D}_r\mathcal{H}$ the space of $\mathbb{Q}$-linear combinations of elements of the form
$\zeta^{\mathfrak{m}}(n_1,\cdots,n_s),\; s\leq r$. Let $gr_r^{\mathfrak{D}}\mathcal{H}=\mathfrak{D}_r\mathcal{H}/\mathfrak{D}_{r+1}\mathcal{H}$ and  $\mathcal{A}=\mathcal{H}/\zeta^{\mathfrak{m}}(2)\mathcal{H}$, then there is an induced depth filtration on $\mathcal{A}$ which is defined by
\[
\mathfrak{D}_r\mathcal{A}=\mathfrak{D}_r\mathcal{H}/\zeta^{\mathfrak{m}}(2)\mathcal{H}\cap\mathfrak{D}_r\mathcal{H}.
\]
Define  $gr_r^{\mathfrak{D}}\mathcal{A}=\mathfrak{D}_r\mathcal{A}/\mathfrak{D}_{r+1}\mathcal{A}$.
Denote by $gr_r^{\mathfrak{D}}\mathcal{A}^{odd}$ the $\mathbb{Q}$-linear subspace of  $gr_r^{\mathfrak{D}}\mathcal{A}$ spanned by the natural image of
\[
\zeta^{\mathfrak{m}}(n_1,\cdots,n_r),\;n_1,\cdots,n_r\geq 3,\; odd
\]
in  $gr_r^{\mathfrak{D}}\mathcal{A}$. For $n_1,\cdots,n_r\geq 3,\; odd$,  $\zeta^{\mathfrak{m}}(n_1,\cdots,n_r)$ is called totally odd motivic multiple zeta value.

From Proposition 10.1 in \cite{depth}, the action of lie algebra $\mathfrak{g}$ on
$\mathfrak{D}_r\mathcal{A}^{odd}/\mathfrak{D}_{r-1}\mathcal{A}$ factors through its abelianization
$\mathfrak{g}^{ab}$, so for $n\geq1$, there is a well-defined derivation
\[\partial_{2n+1}:gr_r^{\mathfrak{D}}\mathcal{A}^{odd}\rightarrow gr_{r-1}^{\mathfrak{D}}\mathcal{A}^{odd}\]
which corresponds to the action of $\widetilde{\sigma}_{2n+1}\in \mathfrak{g}^{ab}$

Here we give the definition of restricted even period polynomials.
\begin{Def}\label{polynomial}
For $N\geq 3$, the restricted even period polynomial of weight $N$ is the polynomial $p(x_1,x_2)$ of degree $N-2$ which satisfies\\
(i) $p(x_1,0)=0$, i.e. p is restricted;\\
(ii) $p(\pm x_1,\pm x_2)=p(x_1,x_2)$, i.e. p is even;\\
(iii) $p(x_1,x_2)+p(x_1-x_2,x_1)-p(x_1-x_2,x_2)=0$.\\
Denote by $\mathbb{P}_N$ the set of even restricted period polynomials of weight $N$.
\end{Def}

Since the action of Lie $G$ on $\mathcal{O}({}_0\Pi_1)$ is compatible with the depth filtration of $\mathcal{O}({}_0\Pi_1)$, we have the following proposition:
\begin{prop}\label{period}
If \[p=\sum_{s+t=N,s,t\geq0}p_{s,t}x_1^{s-1} x_2^{t-1}\]
is a restricted even period polynomial, then we have \\
(i) For a fixed pair $r,s\geq 3$, if at least one of
them is even,
then $p_{r,s}=0$.\\
(ii) By (i), the map
\[
\sum_{s+t=n}p_{s,t}\partial_s\circ\partial_t:gr_r^{\mathfrak{D}}\mathcal{A}^{odd}\rightarrow gr_{r-2}^{\mathfrak{D}}\mathcal{A}^{odd}
\]
is well-defined and we have \[\sum_{s+t=n}p_{s,t}\partial_s\circ\partial_t=0.\]
\end{prop}
\noindent{\bf Proof}:
(i) follows immediately from Definition \ref{polynomial}. (ii) follows from that
\[\sum_{s+t=N}p_{s,t}\{i(\sigma_s),i(\sigma_t)\}\equiv0\]
modulo depth $3$ elements of $\mathbb{L}(e_0,e_1)$, the fact that the action of Lie $G$ on $\mathcal{O}({}_0\Pi_1)$ is compatible with the depth filtration and the even restricted period polynomial is anti-symmetry. While the formula
\[\sum_{s+t=N}p_{s,t}\{i(\sigma_s),i(\sigma_t)\}\equiv0\]
modulo depth $3$ elements of $\mathbb{L}(e_0,e_1)$ follows from Corollary 4.2 in \cite{sch}.  $\hfill\Box$\\

\begin{Def}\label{coation}
We construct a map \[\widetilde{\partial}:gr^{\mathfrak{D}}_r\mathcal{A}^{odd}\rightarrow gr^{\mathfrak{D}}_{1}\mathcal{A}^{odd}\otimes gr^{\mathfrak{D}}_{r-1}\mathcal{A}^{odd}\]
by
\[\widetilde{\partial}(\zeta^{\mathfrak{m}}(n_1,...,n_r))=\sum_{n \,odd  , n\geq3}\zeta^{\mathfrak{m}}(n)\otimes\partial_{n}(\zeta^{\mathfrak{m}}(n_1,...,n_r)).\]
\end{Def}

Essentially the map $\widetilde{\partial}$ is induced from the coaction of $\mathcal{O}(U^{dR})$ on the modulo $\zeta^{\mathfrak{m}}(2)$ version motivic multiple zeta values. The map $\widetilde{\partial}$ is the depth-graded version of the coaction of Lie coalgebra $\mathfrak{g}^{{}^{\vee}}$ on the motivic multiple zeta values (modulo $\zeta^{\mathfrak{m}}(2)$ version) restriced to the totally odd part. See Subsection $3.1$ in \cite{brown} for the coaction of the Lie coalgebra $\mathfrak{g}^{{}^{\vee}}$ on motivic multiple zeta values.

The explicit formula of $\widetilde{\partial}$ was given by the formula (4.4) in \cite{tasaka}. By definition and the motivic action of $\mathfrak{g}$, the operator $\widetilde{\partial}$ preserves the weight.

\begin{Def}\label{vedu}
Denote \[\mathbb{P}^{{}^{\vee}}=\bigoplus_{N\geq3}\mathbb{P}^{{}^{\vee}}_{N},\]
where $\mathbb{P}^{\vee}_N$ is the dual vector space of $\mathbb{P}_N$, elements of $\mathbb{P}^{{}^{\vee}}_N$ are called weight $N$ elements of $\mathbb{P}^{{}^{\vee}}$.

\end{Def}

   Now we construct a map
\[v:gr_1^{\mathfrak{D}}\mathcal{A}^{odd}\otimes_{\mathbb{Q}}gr_1^{\mathfrak{D}}\mathcal{A}^{odd}\rightarrow\mathbb{P}^{\vee}\]
such that $v(\zeta^{\mathfrak{m}}(n_1)\otimes\zeta^{\mathfrak{m}}(n_2))$ is a linear functional whose value at $p=\sum p_{n_1,n_2}x_1^{n_1-1}x_2^{n_2-1}\in\mathbb{P}$ is $p_{n_1,n_2}$. By definition $v$ is surjective.

   Next we construct a map for $r\geq3$
\[D:gr_{1}^{\mathfrak{D}}\mathcal{A}^{odd}\otimes_{\mathbb{Q}}gr_{r-1}^{\mathfrak{D}}\mathcal{A}^{odd}\rightarrow \mathbb{P}^{{}^{\vee}}\otimes gr_{r-2}^{\mathfrak{D}}\mathcal{A}^{odd}\]
such that
\[
D=(v\otimes id_{gr_{r-2}^{\mathfrak{D}}\mathcal{A}^{odd}})\circ( id_{gr_1^{\mathfrak{D}}\mathcal{A}^{odd}}\otimes \widetilde{\partial}).
\]
We will write  $D=(v\otimes id)\circ(id\otimes\widetilde{\partial})$ for convenience.

\section{Nilpotent Lie algebra and unipotent algebraic group}\label{key}

In this section we will give a short introduction to the category equivalence between nilpotent Lie algebras and unipotent algebraic groups. All the results in this section are well-known. The main reference for this section is \cite{cartier}.

Let $\mathfrak{a}$ be a finite dimensional nilpotent Lie algebra over $\mathbb{Q}$, its universal enveloping algebra $\mathcal{U}\mathfrak{a}$ is the tensor algebra $T(\mathfrak{a})$ generated by $\mathfrak{a}$ modulo the two-sided ideal of $T(\mathfrak{a})$ generated by
$x\otimes y-y\otimes x-[x,y],\;x, y\in\mathfrak{a}$,
where $[\,,\,]$ is the Lie bracket of $\mathfrak{a}$.

$\mathcal{U}\mathfrak{a}$ is a $\mathbb{Q}$-algebra with unit element 1. Moreover, we can endow $\mathcal{U}\mathfrak{a}$ with a Hopf algebra structure by the following maps
\[
\Delta:x\mapsto x\otimes1+1\otimes x, S:x\mapsto -x, \epsilon:x\mapsto 0
\]
for all $x\in\mathfrak{a}$.
Denote by $\mathcal{O}=(\mathcal{U}\mathfrak{a})^{{}^{\vee}}$ the compact dual of $\mathcal{U}\mathfrak{a}$. The  elements of $(\mathcal{U}\mathfrak{a})^{{}^{\vee}}$ are the linear functional of $\mathcal{U}\mathfrak{a}$ with finite dimensional support. Then $\mathcal{O}$ is a Hopf algebra with commutative multiplication since the comultiplication of $\mathcal{U}\mathfrak{a}$ is cocommutative.

It's well-known that the above construction gives a category equivalence from the category of finite dimensional nilpotent lie algebras over $\mathbb{Q}$ to the category of Hopf algebra of unipotent algebraic groups over $\mathbb{Q}$. The category equivalence can be generalized to the pro-unipotent case exactly the same way.

We will need the following result in the section \ref{sequences}.
\begin{prop}\label{action}
If $A$ is a pro-unipotent algebraic group over $\mathbb{Q}$ with Lie algebra $\mathfrak{a}=\mathrm{Lie}\, A$, then $\mathfrak{a}$ has a natural action
\[s:\mathfrak{a}\times\mathcal{O}(A)\rightarrow\mathcal{O}(A).\]
The following action of $\mathfrak{a}$
\[r:\mathfrak{a}\times\mathcal{U}\mathfrak{a}\rightarrow\mathcal{U}\mathfrak{a}\]
\[(x,x_1\otimes x_2\otimes\ldots\otimes x_n)\rightarrow x\otimes x_1\otimes x_2\ldots\otimes x_n\]
is well-defined and compatible with the action $s$.
\end{prop}
\noindent{\bf Proof}:
Denote by $\bigtriangleup:\mathcal{O}(A)\rightarrow\mathcal{O}(A)\otimes_{\mathbb{Q}}\mathcal{O}(A)$ and $\epsilon:\mathcal{O}(A)\rightarrow\mathbb{Q}$ the comultiplication and counit map of the Hopf algebra $\mathcal{O}(A)$ respectively. Denote by $I=\mathrm{Ker} \,\epsilon$, we view $\mathbb{Q}$ as the subspace of $\mathcal{O}(A)$ spanned by the unique unit element, then the image of the map $id-\epsilon:\mathcal{O}(A)\rightarrow\mathcal{O}(A)$ is contained in $I$. So we can define $\overline{id-\epsilon}:\mathcal{O}(A)\rightarrow I/I^2$ as the map $id-\epsilon$ modulo $I^2$.

By the classical results of linear algebraic group, we know that $\mathfrak{a}=(I/I^2)^{{}^{\vee}}$ and the map $s$ is determined by the following map
\[((\overline{id-\epsilon})\otimes id)\circ \bigtriangleup :\mathcal{O}(A)\rightarrow(I/I^2)\otimes\mathcal{O}(A).\]
Since the comultiplication $\triangle$ of $\mathcal{O}(A)$ is dual to the multiplication of $\mathcal{U}\mathfrak{a}$, taking dual to the above map we have
\[r:\mathfrak{a}\times\mathcal{U}\mathfrak{a}\rightarrow\mathcal{U}\mathfrak{a},\]
\[(x,x_1\otimes x_2\otimes\ldots\otimes x_n)\rightarrow x\otimes x_1\otimes x_2\ldots\otimes x_n.\]
$r$ is well-defined by the definition of $\mathcal{U}\mathfrak{a}$.$\hfill\Box$\\

Now we consider the simplest unipotent algebraic group $\mathbb{G}_a$ over $\mathbb{Q}$ as an example.

The Lie algebra $\mathfrak{a}$ of $\mathbb{G}_a$ over $\mathbb{Q}$ is one dimensional $\mathbb{Q}$-vector space $\mathbb{Q}x_1$ with trivial Lie bracket. The universal enveloping algebra $\mathcal{U}\mathfrak{a}$ is $\mathcal{U}\mathfrak{a}\cong\mathbb{Q}[x_1]$ with multiplication \[x_1^m\times x_1^n\mapsto x_1^{m+n}\] and comultiplication \[x_1^n\mapsto(x_1\otimes1+1\otimes x_1)^n=\sum_{k=0}^n\dbinom{n}{k}x_1^k\otimes x_1^{n-k}.\] The function ring $\mathcal{O}(\mathbb{G}_a)$ is $\mathcal{O}(\mathbb{G}_a)\cong\mathbb{Q}[x^1]$ with multiplication \[(x^1)^p\times(x^1)^q\mapsto\dbinom{p+q}{p}(x^1)^{p+q}\] and comultiplication
\[(x^1)^n\mapsto\sum_{k=0}^n(x^1)^k\otimes(x^1)^{n-k}.\] The action of $\mathfrak{a}$ on $\mathcal{U}\mathfrak{a}$ is
\[r:\mathfrak{a}\times\mathcal{U}\mathfrak{a}\rightarrow\mathcal{U}\mathfrak{a},\]
\[(x_1,x_1^n)\mapsto x_1^{n+1}.\]
While the action of $\mathfrak{a}$ on $\mathcal{O}(\mathbb{G}_a)$ is
\[s:\mathfrak{a}\times\mathcal{O}(\mathbb{G}_a)\rightarrow\mathcal{O}(\mathbb{G}_a),\]
\[x_1\times(x^1)^q\mapsto(x^1)^{q-1}\]
for $q\geq1$.
\begin{rem}\label{derivation}
Be ware that the ordinary differential of the polynomial ring $\mathbb{Q}[x]$ is \[x^n\mapsto nx^{n-1}.\] But there is no contradiction since the multiplication rule of the ring $\mathcal{O}(\mathbb{G}_a)\cong\mathbb{Q}[x^1]$ is not the ordinary one. The two  differentials that look different will be the same after suitably normalized.
\end{rem}

\section{ Motivic Lie algebra }\label{nonde}

\subsection{A non-degenerated conjecture}
In this subsection we will propose a non-degenerated conjecture in $\mathfrak{d}\mathfrak{g}$. This non-degenerated conjecture will play an important role in proving Theorem \ref{rrr}(ii).

If $V$ is a vector space, denote by $\mathrm{Lie}_n(V)\subseteq V^{\otimes n}$ the degree $n$ part of the free lie algebra generated by $V$. For $n\geq2$, define \[\alpha: \mathbb{P} \otimes\underbrace{\mathfrak{d}\mathfrak{g}_1\otimes\cdots\otimes\mathfrak{d}\mathfrak{g}_1}_{n-2}\rightarrow\mathrm{Lie}_n(\mathfrak{d}\mathfrak{g}_1)\]
 by
\[\alpha:\sum p_{r,s}x_1^{r-1}x_2^{s-1}\otimes \overline{\sigma}_{i_1}\otimes\cdots\otimes\overline{\sigma}_{i_{n-2}}\mapsto  \sum p_{r,s}[\cdots[[\overline{\sigma}_r,\overline{\sigma}_s],\overline{\sigma}_{i_1}],\cdots,\overline{\sigma}_{i_{n-2}}]\]
where $[\,,\,]$ is the formal lie bracket. Denote by $\beta:\mathrm{Lie}_n(\mathfrak{d}\mathfrak{g}_1)\rightarrow \mathfrak{d}\mathfrak{g}_n$  the map that replacing the formal Lie bracket by the induced Ihara bracket. Then we have
\begin{conj}\label{nondege}(non-degenerated conjecture)
For $n\geq2$, the following sequence
\[\mathbb{P} \otimes\underbrace{\mathfrak{d}\mathfrak{g}_1\otimes\cdots\otimes\mathfrak{d}\mathfrak{g}_1}_{n-2}\xrightarrow{\alpha}\mathrm{Lie}_n(\mathfrak{d}\mathfrak{g}_1)\xrightarrow{\beta} \mathfrak{d}\mathfrak{g}_n\]
is exact.

\end{conj}

For $n=2$, this follows from Corollary 4.2 in \cite{sch} since $\mathfrak{dg}_1$ is the space of $\mathbb{Q}$-linear combinations of $(\mathrm{ad}\;e_0)^{2n}e_1,\;n\geq 1$ in $\mathfrak{h}$.

In \cite{depth} Brown introduce linearized double shuffle Lie algebra $\mathfrak{l}\mathfrak{s}\subseteq \mathfrak{h}$ which is also bigraded by weight and depth. By exactly the same method in \cite{racinet} it can be shown that
\[
\mathfrak{d}\mathfrak{g}\subseteq\mathfrak{l}\mathfrak{s}.
\]
By definition $\mathfrak{d}\mathfrak{g}_1=\mathfrak{l}\mathfrak{s}_1$.
For $n=3$, from Theorem 2.6 in \cite{gon} we can deduce that there is an exact sequence
\[
0\rightarrow \mathbb{P}\otimes \mathfrak{ls}_1\rightarrow \mathrm{Lie}_3(\mathfrak{ls}_1)\rightarrow \mathfrak{ls}_3\rightarrow 0
\]
Since $\mathfrak{d}\mathfrak{g}\subseteq\mathfrak{l}\mathfrak{s}$, $\mathfrak{dg}_3=\mathfrak{ls}_3$ and Conjecture \ref{nondege} is true for $n=3$.

The general cases are still open.

\subsection{Multiplication structure on universal enveloping algebra}\label{multiplication}

From the structure of $\mathfrak{h}$, we know that there is a natural isomorphism between the universal enveloping algebra of $\mathfrak{h}$ which is denoted by $\mathcal{U}\mathfrak{h}$ and $\mathbb{Q}\langle e_0,e_1\rangle$:
\[
\lambda:\mathcal{U}\mathfrak{h}\rightarrow\mathbb{Q}\langle e_0,e_1\rangle.
\]
as a $\mathbb{Q}$-vector space. But when transformed to $\mathbb{Q}\langle e_0,e_1\rangle$, the multiplication structure of $\mathcal{U}\mathfrak{h}$ is not the usual concatenation product of $\mathbb{Q}\langle e_0,e_1\rangle$. The new multiplication structure $\circ$ which transformed from $\mathcal{U}\mathfrak{h}$ in some simple cases is calculated explicitly in section 2 and section 6 of Brown \cite{depth}. We review some basic facts of Brown \cite{depth} for convenience.

Denote by $gr^r_{\mathfrak{D}}\mathbb{Q}\langle e_0,e_1\rangle$ the subspace of $\mathbb{Q}\langle e_0,e_1\rangle$ spanned by noncommutative words in $e_0,e_1$ with exactly $r$ occurrences of $e_1$, then we have the following isomorphism (polynomial representation).
\[\rho:gr^r_{\mathfrak{D}}\mathcal{U}\mathfrak{h}\longrightarrow\mathbb{Q}[y_0,...,y_r]\]
\[e_0^{a_0}e_1e_0^{a_1}e_1...e_1e_0^{a_r}\longmapsto y_0^{a_0}y_1^{a_1}...y_r^{a_r}.\]

If $x\in\mathfrak{h}$, then $\lambda(x)$ is the usual Lie polynomial of $x$ in $\mathbb{Q}\langle e_0,e_1\rangle$.

If $x_1,x_2\in\mathfrak{h}$ and $\lambda(x_1),\lambda(x_2)$ belong to $gr^r_{\mathfrak{D}}\mathbb{Q}\langle e_0,e_1\rangle,gr^s_{\mathfrak{D}}\mathbb{Q}\langle e_0,e_1\rangle$ respectively, denote by $\rho(\lambda(x_1))=f(y_0,...,y_r)$ and $\rho(\lambda(x_2))=g(y_0,...,y_s)$, then
\[
\begin{split}
\rho(\lambda(x_1)\circ\lambda(x_2))=&\sum_{i=0}^s f(y_i,y_{i+1},...,y_{i+r})g(y_0,...,y_i,y_{i+r+1},...,y_{r+s})                               \\
                                & +(-1)^{deg f+r}\sum_{i=1}^s f(y_{i+r},...,y_{i+1},y_i)g(y_0,...,y_{i-1},y_{i+r},...,y_{r+s})
\end{split}
\]
\begin{rem}\label{different}
Be aware that the multiplication $\circ$ is not completely the same as the multiplication $\underline{\circ}$
defined by Brown in Definition $2.1$ and Definition $6.2$ in Brown \cite{depth}. The two operations coincide on $\mathfrak{h}\subseteq\mathcal{U}\mathfrak{h}$ (Proposition $2.2$ in Brown \cite{depth}). But in general, $\circ$ satisfies the associative law while $\underline{\circ}$ doesn't.
\end{rem}

\subsection{Totally odd and quasi-uneven}

Since $\mathbb{Q}\langle e^0,e^1\rangle$ is dual to $\mathbb{Q}\langle e_0,e_1\rangle$ and $\mathcal{A}\cong\mathcal{O}(U^{dR})$ is the dual space of  $\mathcal{U}\mathfrak{g}$ and all the depth structures are compatible, we have $gr_r^{\mathfrak{D}}\mathcal{A}\cong gr^r_{\mathfrak{D}}\mathcal{U}\mathfrak{g}$.

  Denote by $\mathnormal{E}$ the linear subspace of $\mathbb{Q}\langle e^0,e^1\rangle$ spanned by the following elements
\[
e^1\underbrace{e^0...e}_{2n_1}{}^0e^1...e^1\underbrace{e^0...e}_{2n_r}{}^0,\;n_1,...,n_r\geq1
\]
Since $\mathcal{A}$ is a quotient space of $\mathbb{Q}\langle e^0,e^1\rangle$, the image of $\mathnormal{E}$ in $\mathcal{A}$ is exactly the modulo $\zeta^{\mathfrak{m}}(2)\mathcal{A}$ version of  totally odd motivic multiple zeta values.

\begin{Def}\label{quasi}
If $x\in gr^r_{\mathfrak{D}}\mathbb{Q}\langle e_0,e_1\rangle$ and the polynomial representation $\rho(x)$ can be written as the form
\[
\rho(x)=g_1(y_0,...,y_r)y_0+g_2(y_1,...,y_r)
\]
where the coefficients of $y_1^{2m_1}y_2^{2m_2}\cdots y_r^{2m_r},m_1,...,m_r\geq1$ in $g_2(y_1,...,y_r)$ are all zero, then we call $x$ quasi-uneven.
\end{Def}

It's clear that all the quasi-uneven elements form a vector space. Furthermore, $\mathnormal{E}$ is the dual space of
$\mathbb{Q}\langle e_0,e_1\rangle$ modulo quasi-uneven elements. So $gr_r^{\mathfrak{D}}\mathcal{A}^{odd}$ is the dual space of  $gr^r_{\mathfrak{D}}\mathcal{U}\mathfrak{g}$ modulo quasi-uneven elements. (Remember that $\mathcal{U}\mathfrak{g}$ is a subalgebra of $\mathcal{U}\mathfrak{h}\cong \mathbb{Q}\langle e_0,e_1\rangle $ with new multiplication.)

Denote by $\gamma:\underbrace{\mathfrak{d}\mathfrak{g}_1\otimes\cdots\otimes\mathfrak{d}\mathfrak{g}_1}_{n}\rightarrow gr^n_{\mathfrak{D}}\mathcal{U}\mathfrak{g}$ the natural map, we propose the following conjecture
\begin{conj}\label{noeven}(vanishing conjecture)
For $n\geq3$, there are no quasi-uneven elements in $\mathrm{Im}\;\gamma$.
\end{conj}

In \cite{depth}, Brown defined uneven elements. But unfortunately uneven elements don't correspond to the totally odd motivic multiple zeta values directly. The space of quasi-uneven elements contains the space of uneven elements. Furthermore, if Conjecture $3$ in \cite{depth} and Conjecture \ref{noeven} are both true, then the space of quasi-uneven elements in $gr_{\mathfrak{D}}^r\mathcal{U}\mathfrak{g}$ is the same as the space of uneven elements in $gr_{\mathfrak{D}}^r\mathcal{U}\mathfrak{g}$.

\begin{rem}
In general, it  should follow that  the uneven and quasi-uneven spaces coincide on the subspace of solutions of the linearized double shuffle equations.
\end{rem}

Denote by $\widetilde{\gamma}:\underbrace{\mathfrak{d}\mathfrak{g}_1\otimes\cdots\otimes\mathfrak{d}\mathfrak{g}_1}_{n}\rightarrow (gr^n_{\mathfrak{D}}\mathcal{U}\mathfrak{g}$  modulo quasi-uneven  elements) the natural quotient map of $\gamma$.

\begin{conj}\label{surconj}(surjective conjecture)
For $n\geq 3$, the map $\widetilde{\gamma}$ is surjective.
\end{conj}

Recall the map $\beta:\mathrm{Lie}_n(\mathfrak{d}\mathfrak{g}_1)\rightarrow\mathfrak{d}\mathfrak{g}_n$, denote by $\mathfrak{d}\mathfrak{g}^{odd}_n$ the $\mathfrak{d}\mathfrak{g}_n$ modulo quasi-uneven elements and $\widetilde{\beta}:\mathrm{Lie}_n(\mathfrak{d}\mathfrak{g}_1)\rightarrow\mathfrak{d}\mathfrak{g}^{odd}_n$ the natural quotient of $\beta$, we have
\begin{Thm}\label{exa}
For $n\geq 3$, assuming Conjecture \ref{nondege}, Conjecture \ref{noeven}, Conjecture \ref{surconj}, there is an exact sequence
\[ \mathbb{P} \otimes\underbrace{\mathfrak{d}\mathfrak{g}_1\otimes\cdots\otimes\mathfrak{d}\mathfrak{g}_1}_{n-2}\xrightarrow{\alpha}\mathrm{Lie}_n(\mathfrak{d}\mathfrak{g}_1)\xrightarrow{\widetilde{\beta}} \mathfrak{d}\mathfrak{g}^{odd}_n\rightarrow 0\]
\end{Thm}

The surjective conjecture is true for $n=3$.
Actually the map $\gamma$ is surjective for $n=3$. From Theorem $2.6$ in \cite{gon}
we have $\mathfrak{ls}_3$ are generated by $\mathfrak{ls}_1$ via the Ihara bracket $\{\;,\;\}$. Since $\mathfrak{dg}_1=\mathfrak{ls}_1$,  $\mathfrak{dg}_3=\mathfrak{ls}_3$ and the map $\gamma$ is surjective for $n=3$.

The three conjectures in this section will play important roles for us to understand Brown's two conjectures (Conjecture \ref{bkb}, Conjecture \ref{un}).

\section{The short exact sequences}\label{sequences}

In this section we will prove Theorem \ref{odd} and give partial results of Conjecture \ref{ttt}.

For $r=2$, 
$v$ is surjective by definition.

From Theorem $2$ in \cite{kaneko} about formal symbols related to double shuffle relation, we can deduce that $gr_2^{\mathfrak{D}}\mathcal{A}^{odd}=gr_2^{\mathfrak{D}}\mathcal{A}$. By taking dual, Im $\partial=$ Ker $v$ follows from Corollary $4.2$ in \cite{sch}, $gr_2^{\mathfrak{D}}\mathcal{A}^{odd}=gr_2^{\mathfrak{D}}\mathcal{A}$  and Proposition \ref{action}.

While  the injectivity of $\widetilde{\partial}$ follows from the fact
that the map with induced Ihara bracket $\mathfrak{d}\mathfrak{g}_1\otimes\mathfrak{d}\mathfrak{g}_1\rightarrow\mathfrak{d}\mathfrak{g}_2:\sigma_r\otimes\sigma_s\mapsto\{\sigma_r,\sigma_s\}$
is suejective.
So Theorem \ref{odd} is proved.

By the results in Section \ref{galois}, we have the following lemma
\begin{lem}\label{contain}
For $r\geq3$, $\mathrm{Im} \;\widetilde{\partial}\subseteq \mathrm{Ker} \;D$, i.e. $D\circ\widetilde{\partial}=0$.
\end{lem}
\noindent{\bf Proof}:
By elementary calculation, we have

\[
\begin{split}
D\circ\widetilde{\partial}&=(v\otimes id)\circ(id\otimes \widetilde{\partial})\circ\widetilde{\partial}\\
                          &=(v\otimes id)\circ\left(\sum_{n_1\;odd,\,\geq3}id\otimes \zeta^{\mathfrak{m}}(n_1)\otimes \partial_{n_1} \right)\circ\left(\sum_{n_2\;odd,\,\geq3}\zeta^{\mathfrak{m}}(n_2)\otimes\partial_{n_2}\right)\\
                          &=(v\otimes id)\circ\left(\sum_{n_1,n_2\;odd,\,\geq3}\zeta^{\mathfrak{m}}(n_2)\otimes\zeta^{\mathfrak{m}}(n_1)\otimes\partial_{n_1}\circ\partial_{n_2}\right)\\
                          &=\sum_{n_1,n_2\;odd,\,\geq3}v(\zeta^{\mathfrak{m}}(n_2)\otimes\zeta^{\mathfrak{m}}(n_1))\otimes \partial_{n_1}\circ\partial_{n_2} \\
                          &=0
\end{split}
\]
The last equality is derived from the definition of $v$ and Proposition \ref{period}.$\hfill\Box$\\

Since $\mathcal{A}\cong\mathcal{O}(U^{dR})$ is the dual space of  $\mathcal{U}\mathfrak{g}$ and the depth filtration on $\mathcal{A}$ and $\mathfrak{g}$ is compatible, so we have $gr_r^{\mathfrak{D}}\mathcal{A}$ is the dual space of $gr^r_{\mathfrak{D}}\mathcal{U}\mathfrak{g}$. Since $gr_r^{\mathfrak{D}}\mathcal{A}^{odd}$ is the totally odd part of $gr_r^{\mathfrak{D}}\mathcal{A}$, we have $gr_r^{\mathfrak{D}}\mathcal{A}^{odd}$ is the dual space of  $gr_{\mathfrak{D}}^r\mathcal{U}\mathfrak{g}$ modulo quasi-uneven elements.

\begin{prop}\label{equal}
If Conjecture \ref{nondege} (non-degenarated conjecture), Conjecture \ref{noeven}\\(vanishing conjecture) and Conjecture \ref{surconj}(surjective conjecture) are true for all $r\geq3$, then for each $r\geq 3$ there is an exact sequence
\[
 0\rightarrow gr^{\mathfrak{D}}_r\mathcal{A}^{odd}\xrightarrow{\widetilde{\partial}} gr^{\mathfrak{D}}_{1}\mathcal{A}^{odd}\otimes gr^{\mathfrak{D}}_{r-1}\mathcal{A}^{odd}\xrightarrow{D} \mathbb{P}^{{}^{\vee}}\otimes gr^{\mathfrak{D}}_{r-2}\mathcal{A}^{odd}\rightarrow 0.
\]
 .
\end{prop}
\noindent{\bf Proof}:
Since $gr_r^{\mathfrak{D}}\mathcal{A}^{odd}$ is the totally odd
 part of $gr_r^{\mathfrak{D}}\mathcal{A}$, by taking dual, the injectivity of $\widetilde{\partial}$ follows from the assumption (Conjecture \ref{surconj})  that $\widetilde{\gamma}$ is surjective and Proposition \ref{action}.

 From the surjective conjecture and vanishing conjecture, we know the dimension of the weight $N$ part of
 $(gr_{\mathfrak{D}}^r\mathcal{U}\mathfrak{g})$ modulo quasi-uneven elements is equal to the dimension of the weight $N$ part of the natural image $\underbrace{\mathfrak{d}\mathfrak{g}_1\otimes\cdots\otimes\mathfrak{d}\mathfrak{g}_1}_{r}$ in  $gr^r_{\mathfrak{D}}\mathcal{U}\mathfrak{g}$. So from the non-degenerated conjecture, $\mathrm{Im}\;\widetilde{\partial}=\mathrm{Ker} \;D$ .
 
 We regard $\mathbb{P}$ as a subspace of $\mathfrak{d}\mathfrak{g}_1\otimes\mathfrak{d}\mathfrak{g}_1$ via the following map
\[t:\mathbb{P}\rightarrow\mathfrak{d}\mathfrak{g}_1\otimes\mathfrak{d}\mathfrak{g}_1\]
\[\sum_{\substack{r+s=N\\r,s \,odd,\,\geq3}}p_{r,s}x^{r-1}y^{s-1}\mapsto\sum_{\substack{r+s=N\\r,s \,odd,\,\geq3}}p_{r,s}\overline{\sigma}_r\otimes\overline{\sigma}_s\]
Then by the work of Goncharov \cite{goncharov}, we know that
\[\mathbb{P}\otimes\mathfrak{d}\mathfrak{g}_1\cap\mathfrak{d}\mathfrak{g}_1\otimes\mathbb{P}=\{0\}\]
in $\mathfrak{d}\mathfrak{g}_1\otimes\mathfrak{d}\mathfrak{g}_1\otimes\mathfrak{d}\mathfrak{g}_1$.
 
 Denote by $\mathcal{U}\mathfrak{dg}_1$ the universal enveloping algebra of the Lie subalgebra of $\mathfrak{dg}$ generated by $\mathfrak{dg}_1$.  Assuming   Conjecture \ref{nondege}, Conjecture \ref{noeven} and Conjecture \ref{surconj} for all $r\geq 3$, then the dual map $D^{{}^{\vee}}$ is essentially the map 
 \[
 D^{{}^{\vee}}: \mathbb{P}\otimes gr^{r-2}_{\mathfrak{D}}\mathcal{U}\mathfrak{dg}_1\to \mathfrak{dg}_1\otimes gr_{\mathfrak{D}}^{r-1}\mathcal{U}\mathfrak{dg}_1,
 \]
 \[
 \sum p_{r,s} x^{r-1} y^{s-1}\otimes b\mapsto \sum p_{r,s} \overline {\sigma}_r\otimes (\overline{\sigma}_s\otimes b),
 \]
 and the dual map $\widetilde{\partial}^{{}^{\vee}}$ is essentially the map 
 \[
 \widetilde{\partial}^{{}^{\vee}}: \mathfrak{dg}_1\otimes gr_{\mathfrak{D}}^{r-1}\mathcal{U}\mathfrak{dg}_1\to gr_{\mathfrak{D}}^{r}\mathcal{U}\mathfrak{dg}_1,
 \]
 \[
 \overline{\sigma}_r\otimes c\mapsto (\overline{\sigma}_r\otimes c).
 \]
 
 Define $$K: \mathbb{P}\otimes gr_{\mathfrak{D}}^{r-2}\mathcal{U}\mathfrak{dg}_1\to \mathfrak{dg}_1\otimes \mathfrak{dg}_1\otimes gr_{\mathfrak{D}}^{r-2}\mathcal{U}\mathfrak{dg}_1$$ as 
 $$K(\sum_{\substack{r+s=N\\r,s \,odd,\,\geq3}}p_{r,s}x^{r-1}y^{s-1}\otimes d)=\sum_{\substack{r+s=N\\r,s \,odd,\,\geq3}}p_{r,s}\overline{\sigma}_r\otimes\overline{\sigma}_s\otimes d.$$
 It's clear that 
 \[
 D^{{}^{\vee}}=(id|_{\mathfrak{dg}_1}\otimes \widetilde{\partial}^{{}^{\vee}})\circ K.
 \]
 
  Since $K$ is injective, from the equality $\mathrm{Ker}\; \widetilde{\partial}=\mathrm{Im}\; D^{{}^{\vee}}$ in depth $r-1$ we have 
  \[
  \begin{split}
  \mathrm{Ker}\; D^{{}^{\vee}}&= \mathrm{Ker}\;(id|_{\mathfrak{dg}_1}\otimes \widetilde {\partial}^{{}^{\vee}})\cap \mathrm{Im}\; K\\
                                                  &=\mathrm{Im}\; (id|_{\mathfrak{dg}_1}\otimes D^{{}^{\vee}})\cap \mathrm{Im}\;K\\
                                                  &\subseteq (id|_{\mathfrak{dg}_1}\otimes D^{{}^{\vee}}) \left ((\mathbb{P}\otimes \mathfrak{dg}_1\cap \mathfrak{dg}_1\otimes \mathbb{P})\otimes gr_{\mathfrak{D}}^{r-3}\mathcal{U}\mathfrak{dg}_1   \right).
  \end{split}
  \]
So from $\mathbb{P}\otimes\mathfrak{d}\mathfrak{g}_1\cap\mathfrak{d}\mathfrak{g}_1\otimes\mathbb{P}=\{0\}$, we deduce that $D^{{}^{\vee}}$ is injective in depth $r$ and $D$ is surjective in depth $r$.  $\hfill\Box$\\

Now  we discuss the map $D$ using explicit matrix calculation, we first introduce some notation from Tasaka\cite{tasaka}.
Denote by
\[
S_{N,r}=\{(n_1,...,n_r)\in\mathbb{Z}^r\mid n_1+...+n_r=N,n_1,...,n_r\geq3:odd\}.
\]
We write $\overrightarrow{m}=(m_1,...,m_r)$ for short, while \[\mathbf{Vect}_{N,r}=\{(a_{n_1,...,n_r})_{\overrightarrow{n}\in S_{N,r}}\mid a_{n_1,...,n_r}\in \mathbb{Q}\}.\]

For a matrix $P=\left(p\dbinom{m_1,...,m_r}{n_1,...,n_r}\right)_{\substack {\overrightarrow{m}\in S_{N,r}\\\overrightarrow{n}\in S_{N,r}}}$, the action of $P$ on $a=(a_{m_1,...,m_r})_{\overrightarrow{m}\in S_{N,r}}$ means
\[
aP=\left(\sum_{\overrightarrow{m}\in S_{N,r}}a_{m_1,...,m_r}p\dbinom{m_1,...,m_r}{n_1,...,n_r}\right)_{\overrightarrow{n}\in S_{N,r}}
\]
Denote by $\mathbb{P}_{N,r}$ the $\mathbb{Q}$-vector space spanned by the set
\[
\{x_1^{n_1-1}\cdots x_r^{n_r-1}\mid(n_1,...,n_r)\in S_{N,r}\}.
\]
Obviously there is an isomorphism
\[\pi_1:\mathbb{P}_{N,r}\longrightarrow\mathbf{Vect}_{N,r}\]
\[\sum_{\overrightarrow{n}\in S_{N,r}}a_{n_1,...,n_r}x^{n_1-1}_1\cdots x^{n_r-1}_r\longmapsto(a_{n_1,...,n_r})_{\overrightarrow{n}\in S_{N,r}}.\]
Denote by
\[
\mathbf{W}_{N,r}=\{p\in\mathbb{P}_{N,r}\mid p(x_1,...,x_r)=p(x_2-x_1,x_2,x_3,...,x_r)-p(x_2-x_1,x_1,x_3,...,x_r)\}.
\]
For $s\geq 1$, define $\delta\dbinom{a_1,...,a_s}{b_1,...,b_s}=1$ if $(a_1,\cdots,a_s)=(b_1,\cdots,b_s)$, $\delta\dbinom{a_1,\cdots,a_s}{b_1,\cdots, b_s}=0$ if $(a_1,\cdots,a_s)\neq(b_1,\cdots,b_s)$.

From the formula $(4.4)$ in \cite{tasaka}, we have
\[
\widetilde{\partial}(\zeta^{\mathfrak{m}}(n_1,...,n_r))=\sum_{\overrightarrow{m}\in S_{N,r}}e\dbinom{m_1,...,m_r}{n_1,...,n_r}\zeta^{\mathfrak{m}}(m_1)\otimes\zeta^{\mathfrak{m}}(m_2,...,m_r).
\]
Where
\[
e\dbinom{m_1,...,m_r}{n_1,...,n_r}=\delta\dbinom{m_1,...,m_r}{n_1,...,n_r}+\sum_{i=1}^{r-1}\delta\dbinom{m_2,...,m_i,m_{i+2},...,m_r}{n_1,...,n_{i-1},n_{i+2},...,n_r}b^{m_1}_{n_i,n_{i+1}},
\]
and the $b^m_{n,n'}$ are defined by
\[
b^m_{n,n'}=(-1)^n\dbinom{m-1}{n-1}+(-1)^{n'-m}\dbinom{m-1}{n'-1}.
\]
The above formula of $\widetilde{\partial}$ can be deduced from Proposition $2.2$ and the calculation of Subsection $6.1$ in \cite{bro1}.

The matrix $E^{(r-i)}_{N,r},\,i=1,...,r-2$ are defined by
\[
E^{(r-i)}_{N,r}=\left(\delta\dbinom{m_1,...,m_i}{n_1,...,n_i}e\dbinom{m_{i+1},...,m_r}{n_{i+1},...,n_r}\right)_{\substack {\overrightarrow{m}\in S_{N,r}\\\overrightarrow{n}\in S_{N,r}}}.
\]
Denote by
\[
C_{N,r}=E^{(2)}_{N,r}\cdot E^{(3)}_{N,r}\cdots E^{(r-1)}_{N,r}\cdot E_{N,r}
\]

In \cite{tasaka}, Tasaka proposed the following conjecture
\begin{conj}\label{general}
For $r\geq 3$, the map
\[\eta:\mathrm{W}_{N,r}\rightarrow \mathrm{V}_{N,r}
\]
\[
p(x_1,...,x_r)\mapsto (\pi_1(p(x_1,...,x_r))(E_{N,r}-I_{N,r})
\]
is injective.
\end{conj}

Now we will prove Conjecture \ref{general} in case $r=3$.
\begin{prop}\label{spec}
The map
\[\eta:\mathrm{W}_{N,3}\rightarrow \mathrm{V}_{N,3}
\]
\[
p(x_1,x_2,x_3)\rightarrow (\pi_1(p(x_1,x_2,x_3))(E_{N,3}-I_{N,3})
\]
is injective.
\end{prop}
\noindent{\bf Proof}:
If $p\in \mathrm{W}_{N,3}$ satisfies $\eta(p)=0$, denote
\[
p=\sum_{\overrightarrow{m}\in S_{N,3}}a_{m_1,m_2,m_3}x_1^{m_1-1}x_2^{m_2-1}x_3^{m_3-1},
\]
then we have
\[
a_{m_1,m_2,m_3}+a_{m_2,m_1,m_3}=0
\]
\[
\sum_{\overrightarrow{m}\in S_{N,3}}a_{m_1,m_2,m_3}\left(\delta\binom{m_3}{n_3}b^{m_1}_{n_1,n_2}+\delta\binom{m_2}{n_1}b^{m_1}_{n_2,n_3}\right)=0
\]
\[
\sum_{\overrightarrow{m}\in S_{N,3}}a_{m_1,m_2,m_3}\left(\delta\dbinom{m_1,m_2,m_3}{n_1,n_2,n_3}+\delta\binom{m_3}{n_3}b^{m_1}_{n_1,n_2}\right)=0.
\]
So
\[
a_{n_1,n_2,n_3}=\sum_{\overrightarrow{m}\in S_{N,3}}a_{m_1,m_2,m_3}\delta\dbinom{m_1,m_2,m_3}{n_1,n_2,n_3}=\sum_{\overrightarrow{m}\in S_{N,3}}a_{m_1,m_2,m_3}\delta\binom{m_2}{n_1}b^{m_1}_{n_2,n_3}.
\]

Since $b^{m_1}_{n_2,n_3}+b^{m_1}_{n_3,n_2}=0$, so $a_{n_1,n_2,n_3}+a_{n_1,n_3,n_2}=0$. We identify $\mathbb{P}\otimes\mathfrak{d}\mathfrak{g}_1, \mathfrak{d}\mathfrak{g}_1\otimes\mathfrak{d}\mathfrak{g}_1\otimes\mathfrak{d}\mathfrak{g}_1$ with $W_{N,3},\;V_{N,3}$ respectively in the natural way, then we have
\[
p\in\mathbb{P}\otimes\mathfrak{d}\mathfrak{g}_1\cap\wedge^3\mathfrak{d}\mathfrak{g}_1.
\]

By the main theorem of Goncharov \cite{gon}, $\mathbb{P}\otimes\mathfrak{d}\mathfrak{g}_1\cap\wedge^3\mathfrak{d}\mathfrak{g}_1=\{0\}$, so $p=0$.   $\hfill\Box$\\

Now we have
\begin{prop}\label{sur}
For $r=3$, the map $D$ is surjective. For $r>3$, if Conjecture \ref{general} is true, then $D$ is surjective.
\end{prop}
\noindent{\bf Proof}:
Since the explicit formula of the map $D$ is
\[
D(\zeta^{\mathfrak{m}}(n_1)\otimes\zeta^{\mathfrak{m}}(n_2,...,n_{r}))=\sum_{\overrightarrow{m}\in S_{N,r}} e\dbinom{m_2,...,m_{r}}{n_2,...,n_{r}}v(\zeta^{\mathfrak{m}}(n_{1}),\zeta^{\mathfrak{m}}(m_2))\otimes \zeta^{\mathfrak{m}}(m_3,\ldots,m_{r}).
\]
To prove that $D$ is surjective, it suffices to prove that the following map $\xi$ is injective. (Be careful that $\xi$ is not the dual map of $D$, but $D$ is a quotient of the dual map of $\xi$ in some sense)

\[\xi:\mathbf{W}_{N,r}\longrightarrow \mathbf{Vect}_{N,r}\]
\[p(x_1,...,x_r)\longmapsto (\pi_1(p(x_1,...,x_r))E^{(r-1)}_{N,r}.\]
i.e., we need to prove that $\pi_1 (\mathbf{W}_{N,r})\cap \mathrm{ker}\,\mathrm E^{(r-1)}_{N,r}={0}$.

Set $\pi_1(p)=(a_{n_1,...,n_r})_{\overrightarrow{n}\in S_{N,r}}$ for $p\in \mathbf{W}_{N,r}$, the assumption $p\in\mathbf{W}_{N,r}$ implies

\[
\sum_{\overrightarrow{m}\in S_{N,r}}\left(\delta\dbinom{m_1,...,m_r}{n_1,...,n_r}+\delta\dbinom{m_3,...,m_r}{n_3,...,n_r}b^{m_1}_{n_1,n_2}\right)a_{m_1,...,m_r}=0.   \tag{1}
\]

If $(\pi_1(p))E^{(r-1)}_{N,r}=0$, we will have
\[
\sum_{\overrightarrow{m}\in S_{N,r}}\left(\delta\dbinom{m_1,...,m_r}{n_1,...,n_r}+\sum_{i=2}^{r-1}\delta\dbinom{m_1,m_3,...,m_i,m_{i+2},...,m_r}{n_1,n_2,...,n_{i-1},n_{i+2},...,n_r}b^{m_2}_{n_i,n_{i+1}} \right)a_{m_1,...,m_r}=0.  \tag{2}
\]

Since $(a_{n_1,...,n_r})_{\overrightarrow{n}\in S_{N,r}}\in \pi_1(W_{N,r})$ we have
\[
a_{m_1,m_2,m_3,...,m_r}+a_{m_2,m_1,m_3,...,m_r}=0. \tag{3}
\]

From the formulas $(1)$ and  $(3)$, we have
\[
\sum_{\overrightarrow{m}\in S_{N,r}}\left(
\delta\dbinom{m_1,...,m_r}{n_1,...,n_r}+\delta\dbinom{m_3,...,m_r}{n_3,...,n_r}b^{m_2}_{n_2,n_1}\right)a_{m_1,...,m_r}=0.  \tag{4}
\]

From the formulas $(2)$ and  $(4)$ and the fact that $b^{m_2}_{n_2,n_1}+b^{m_2}_{n_1,n_2}=0$, we have
\[
(\pi_1(p))(E_{N,r}-I)=0.
\]
So by Proposition \ref{spec} and Conjecture \ref{general}, the proposition is proved.         $\hfill\Box$\\

\begin{rem}\label{geometric}
Applying the method in \cite{tasaka} we can prove the following result in the proof of Proposition \ref{sur} :
\[
\mathrm{Im} \,\xi\subseteq\mathrm{Ker} \,E^{(r-1)}_{N,r}E_{N,r}.
\]
Using the above result we immediately obtain a second proof of Lemma \ref{contain} based on explicit calculation. So the maps we have established provide a geometric explanation of the results in \cite{tasaka}.
\end{rem}

\section{ Further conjecture}\label{app}

Now we consider the general cases. Exactly the same way as the totally odd case, we can construct the following map
\[\widetilde{\partial}_g:gr^{\mathfrak{D}}_r\mathcal{A}\rightarrow gr^{\mathfrak{D}}_{1}\mathcal{A}^{odd}\otimes gr^{\mathfrak{D}}_{r-1}\mathcal{A}\]
\[D_g:gr_{1}^{\mathfrak{D}}\mathcal{A}^{odd}\otimes gr_{r-1}^{\mathfrak{D}}\mathcal{A}\rightarrow \mathbb{P}^{{}^{\vee}}\otimes gr_{r-2}^{\mathfrak{D}}\mathcal{A}.\]
We have $\widetilde{\partial}_g\mid_{gr_1^{\mathfrak{D}}\mathcal{A}^{odd}}=\widetilde{\partial}$,
$D_g\mid_{gr_1^{\mathfrak{D}}\mathcal{A}^{odd}\otimes gr_{r-1}^{\mathfrak{D}}\mathcal{A}^{odd}}=D$.
This means that $\widetilde{\partial}_g$ and $D_g$ are the general cases of the map $\widetilde{\partial}$ and $D$ respectively.

\begin{prop}\label{half}
(i)For $r=3$, there is an exact sequence
 \[0\rightarrow gr^{\mathfrak{D}}_3\mathcal{A}\xrightarrow{\widetilde{\partial}_g} gr^{\mathfrak{D}}_{1}\mathcal{A}^{odd}\otimes gr^{\mathfrak{D}}_{2}\mathcal{A}\xrightarrow{D_g} \mathbb{P}^{{}^{\vee}}\otimes gr^{\mathfrak{D}}_{1}\mathcal{A}\rightarrow 0\]
 (ii)For $r>3$, $\mathrm{Im}\;\widetilde {\partial}_g\subseteq\mathrm{Ker} \;D_g$. If Conjecture \ref{nondege} is true for all $r\geq 2$, then there is an exact sequence
 \[
 gr^{\mathfrak{D}}_r\mathcal{A}\xrightarrow{\widetilde{\partial}_g} gr^{\mathfrak{D}}_{1}\mathcal{A}^{odd}\otimes gr^{\mathfrak{D}}_{r-1}\mathcal{A}\xrightarrow{D_g} \mathbb{P}^{{}^{\vee}}\otimes gr^{\mathfrak{D}}_{r-2}\mathcal{A}\rightarrow 0\]

\end{prop}
\noindent{\bf Proof}:
From Theorem $2.6$ (b) in \cite{gon}, we can deduce that there is an exact sequence
\[
0\rightarrow \mathbb{P}\otimes \mathfrak{ls}_1\rightarrow \mathrm{Lie}_3(\mathfrak{ls}_1)\rightarrow \mathfrak{ls}_3\rightarrow 0
\]
Since $\mathfrak{ls}_1=\mathfrak{dg}_1$ and $\mathfrak{dg}\subseteq\mathfrak{ls}$, we have $\mathfrak{dg}_3=\mathfrak{ls}_3$ and the following exact sequence
\[
0\rightarrow \mathbb{P}\otimes \mathfrak{dg}_1\rightarrow \mathrm{Lie}_3(\mathfrak{dg}_1)\rightarrow \mathfrak{dg}_3\rightarrow 0
\]

Since $\mathcal{A}\cong\mathcal{O}(\mathrm{U}^{dR})$,
by taking dual, $(i)$ follows immediately from Proposition \ref{action} and the above exact sequence.

For $r>3$, assuming Conjecture \ref{nondege},  by taking dual it suffices to prove the following sequence
\[
0\rightarrow\mathbb{P}\otimes gr_{\mathfrak{D}}^{r-2}\mathcal{U}\mathfrak{g}\xrightarrow{D_g^{{}^{\vee}}}\mathfrak{d}\mathfrak{g}_1\otimes gr_{\mathfrak{D}}^{r-1}\mathcal{U}\mathfrak{g}\xrightarrow{\widetilde{\partial}_g^{{}^{\vee}}} gr_{\mathfrak{D}}^r\mathcal{U}\mathfrak{g}
\]
is exact.
    By Proposition \ref{action}, $\widetilde{\partial}_g^{{}^{\vee}}$ is induced from the following map
\[
\mathfrak{g}\otimes\mathcal{U}\mathfrak{g}\rightarrow\mathcal{U}\mathfrak{g},
\]
\[
a_1\otimes(a_2\otimes\cdots\otimes a_r)\mapsto a_1\otimes a_2\otimes\cdots\otimes a_r,
\]
and $D_g^{{}^{\vee}}$ is induced from the following map
\[
\mathbb{P}\otimes \mathcal{U}\mathfrak{g}\rightarrow \mathfrak{g}\otimes\mathfrak{g}\otimes \mathcal{U}\mathfrak{g}\rightarrow \mathfrak{g}\otimes \mathcal{U}\mathfrak{g},
\]
\[
\sum_{\substack{r,s\geq 3,\\odd}}p_{r,s}x_1^{r-1}x_2^{s-1}\otimes b \mapsto \sum_{\substack{r,s\geq 3,\\odd}}p_{r,s}\sigma_r\otimes\sigma_s\otimes b \mapsto\sum_{\substack{r,s\geq 3,\\odd}}p_{r,s}\sigma_r\otimes (\sigma_s\otimes b).
\]

Denote by $\mathfrak{d}\mathfrak{g}_i$ the depth $i$ elements of $\mathfrak{d}\mathfrak{g}$.
Let $a\in \mathfrak{d}\mathfrak{g}_1\otimes gr_{\mathfrak{D}}^{r-1}\mathcal{U}\mathfrak{g}$. If $a\in \mathrm{Ker}\;\widetilde{\partial}_g^{{}^{\vee}}$, then by the definition of the depth filtration on $\mathcal{U}\mathfrak{g}$, a lift of $a$ in $\mathfrak{g}\otimes \mathcal{U}\mathfrak{g}$ should be the form
\[
\widetilde{a}=\sum\{\cdots\{\widetilde{c}_{k_1},\widetilde{c}_{k_2}\},\cdots,\widetilde{c}_{k_i}\}\otimes \widetilde{d}_{k_{i+1}},
\]
			where $\widetilde{c}_{k_i}\in \mathfrak{D}^{m_{k_i}}\mathfrak{g}$, $\widetilde{d}_{k_{i+1}}\in \mathfrak{D}^{r-m_{k_1}-\cdots-m_{k_i}}\mathcal{U}\mathfrak{g}$ and
			\[
			\{\cdots\{\widetilde{c}_{k_1},\widetilde{c}_{k_2}\},\cdots,\widetilde{c}_{k_i}\}\in \mathfrak{D}^{m_{k_1}+\cdots+m_{k_i}+1}\mathfrak{g}.\]
			Be aware that $\mathfrak{g}$ is viewed as a Lie subalgebra of $\mathcal{U}\mathfrak{g}$ and $\mathfrak{g}\otimes \mathcal{U}\mathfrak{g}$ is viewed as a subspace of $\mathcal{U}\mathfrak{g}\otimes \mathcal{U}\mathfrak{g}$.

Denote by $c_{k_i}$ the image of $\widetilde{c}_{k_i}$ in $\mathfrak{d}\mathfrak{g}_{m_{k_i}}$ and $d_{k_{i+1}}$ the image of $\widetilde{d}_{k_{i+1}}$ in $$gr_{\mathfrak{D}}^{r-m_{k_1}-\cdots-m_{k_i}}\mathcal{U}\mathfrak{g},$$ then $c_{k_i}$  must lie in the image of $\beta: \mathrm{Lie}_{m_{k_i}}(\mathfrak{dg}_1)\rightarrow \mathfrak{d}\mathfrak{g}_{m_{k_i}}$. Otherwise
\[
a=\sum \mathrm{Lie\; polynomial\; of}\;\{\cdots\{{c}_{k_1},{c}_{k_2}\},\cdots,{c}_{k_i}\}\otimes{d}_{k_{i+1}}
\]
can't lie in $\mathfrak{d}\mathfrak{g}_1\otimes gr_{\mathfrak{D}}^{r-1}\mathcal{U}\mathfrak{g}$.
So if $a \in \mathrm{Ker}\;\widetilde{\partial}_g^{{}^{\vee}}$, assuming Conjecture \ref{nondege}, we have $a\in\mathrm{Im}\;D_g^{{}^{\vee}}$.

We view $\mathbb{P}$ as a subspace of $\mathfrak{dg}_1\otimes \mathfrak{dg}_1$ the same way as in the proof of Proposition \ref{equal}.
By the work of Goncharov \cite{goncharov}, 
\[\mathbb{P}\otimes\mathfrak{d}\mathfrak{g}_1\cap\mathfrak{d}\mathfrak{g}_1\otimes\mathbb{P}=\{0\}\]
in $\mathfrak{d}\mathfrak{g}_1\otimes\mathfrak{d}\mathfrak{g}_1\otimes\mathfrak{d}\mathfrak{g}_1$.

Define $$J:\mathbb{P}\otimes gr_{\mathfrak{D}}^{r-2}\mathcal{U}\mathfrak{g}\to \mathfrak{dg}_1\otimes \mathfrak{dg}_1\otimes gr_{\mathfrak{D}}^{r-2}\mathcal{U}\mathfrak{g}$$ as 
$$   J(\sum_{\substack{r+s=N\\r,s \,odd,\,\geq3}}p_{r,s}x^{r-1}y^{s-1}\otimes b)=\sum_{\substack{r+s=N\\r,s \,odd,\,\geq3}}p_{r,s}\overline{\sigma}_r\otimes\overline{\sigma}_s\otimes b  .$$
It's clear that 
\[
D_g^{{}^{\vee}}=(id|_{\mathfrak{dg}_1}\otimes \widetilde{\partial}_g^{{}^{\vee}})\circ J.
\]

Since $J$ is injective, from the equality $\mathrm{Ker}\; \widetilde{\partial}_g^{{}^{\vee}}=\mathrm{Im}\; D_g^{{}^{\vee}}$ in depth $r-1$ we have 
\[
\begin{split}
\mathrm{Ker}\; D_g^{{}^{\vee}}&= \mathrm{Ker}\;(id|_{\mathfrak{dg}_1}\otimes \widetilde{\partial}_g^{{}^{\vee}})\cap \mathrm{Im}\; J\\
                                                     &=\mathrm{Im}\;(id|_{\mathfrak{dg}_1}\otimes D_g^{{}^{\vee}}   )\cap \mathrm{Im}\; J\\
                                                     &\subseteq (id|_{\mathfrak{dg}_1}\otimes D_g^{{}^{\vee}})\left ((\mathbb{P}\otimes\mathfrak{d}\mathfrak{g}_1\cap\mathfrak{d}\mathfrak{g}_1\otimes\mathbb{P})\otimes gr_{\mathfrak{D}}^{r-3}\mathcal{U}\mathfrak{g}\right).
\end{split}
\]
Since 
\[\mathbb{P}\otimes\mathfrak{d}\mathfrak{g}_1\cap\mathfrak{d}\mathfrak{g}_1\otimes\mathbb{P}=\{0\},\]
$D_g^{{}^{\vee}}$ is injective in depth $r$, and $D_g$ is surjective in depth $r$.   $\hfill\Box$\\

Inspired by Theorem \ref{half}, we propose the following conjecture.
\begin{conj}\label{long}
For $r\geq4$, there is an exact sequence
\[
0\rightarrow\mathbb{P}^{{}^{\vee}}\otimes gr_{r-4}^{\mathfrak{D}}\mathcal{A}\rightarrow gr_r^{\mathfrak{D}}\mathcal{A}\rightarrow gr_1^{\mathfrak{D}}\mathcal{A}^{odd}\otimes gr_{r-1}^{\mathfrak{D}}\mathcal{A}\rightarrow \mathbb{P}^{{}^{\vee}}\otimes gr_{r-2}^{\mathfrak{D}}\mathcal{A}\rightarrow0.
\]

\end{conj}
The first map is not clear enough at present. The exceptional elements which are constructed by Brown in Section $8$\cite{depth} should be useful to understand the first map if they are motivic.

In \cite{zeta3},  Brown constructed a map $\mathfrak{c}:\mathbb{P}\rightarrow \mathfrak{d}\mathfrak{g}_4$.
The map $\mathfrak{c}$ should be enough to give the first map. Unfortunately the explicit formula of $\mathfrak{c}$ is  too complicated. The author doesn't know how to deduce the explicit formula of the first map from $\mathfrak{c}$.

We have the following theorem
\begin{Thm}\label{big}
The motivic Broadhurst-Kreimer conjecture follows from the long exact sequence conjecture.
\end{Thm}
\noindent{\bf Proof}:
For $r=1$ it's clear that
\[
\sum_{N>0}\mathrm{dim}_{\mathbb{Q}}(gr_1^{\mathfrak{D}}\mathcal{H}_N)x^N=\mathbb{O}(x)+\mathbb{E}(x).
\]

From Corollary $8$ (parity results) in \cite{ihara}, since the proof of parity results only use double shuffle relations, we know that $gr_2^{\mathfrak{D}}\mathcal{A}_N=0$ if $N$ is odd. From Theorem $1$, $2$ in \cite{kaneko}, we can deduce that $gr_2^{\mathfrak{D}}\mathcal{A}_N=gr_2^{\mathfrak{D}}\mathcal{A}_N^{odd}$  for $N$ even by the same reason. The isomorphism \[\mathcal{H}\cong\mathcal{O}(\mathrm{U}^{dR})\otimes\mathbb{Q}[\zeta^{\mathfrak{m}}(2)]\]
implies 
\[
gr_2^{\mathfrak{D}}\mathcal{H}\cong gr_2^{\mathfrak{D}}\mathcal{A}\oplus \mathbb{Q}[\zeta^{\mathfrak{m}}(2)]\otimes gr_1^{\mathfrak{D}}\mathcal{A}.
\]
So
\[
\begin{split}
\sum_{N>0}\mathrm{dim}_{\mathbb{Q}}(gr_2^{\mathfrak{D}}\mathcal{H}_N)x^N&=\sum_{N>0}\mathrm{dim}_{\mathbb{Q}}(gr_2^{\mathfrak{D}}\mathcal{A}_N)x^N+\sum_{N>0}\mathrm{dim}_{\mathbb{Q}}(gr_1^{\mathfrak{D}}\mathcal{A}_N)x^N\cdot\frac{1}{1-x^2}                                                       \\
                                                                                           &=\sum_{N>0}\mathrm{dim}_{\mathbb{Q}}(gr_2^{\mathfrak{D}}\mathcal{A}_N) x^N+\mathbb{O}(x)\cdot \mathbb{E}(x)
\end{split}
\]
From Theorem \ref{odd}, we have
\[
\sum_{N>0}\mathrm{dim}_{\mathbb{Q}} (gr_2^{\mathfrak{D}}\mathcal{A}_N^{odd})x^N-\left(\sum_{N>0}\mathrm{dim}_{\mathbb{Q}}(gr_1^{\mathfrak{D}}\mathcal{A}_N)x^N\right)^2+\sum_{N>0}\mathrm{dim}_{\mathbb{Q}}(\mathbb{P}^{{}^{\vee}}_N )x^N=0.
\]
So
\[
\sum_{N>0}\mathrm{dim}_{\mathbb{Q}} (gr_2^{\mathfrak{D}}\mathcal{A}_N)x^N=\mathbb{O}(x)^2-\mathbb{S}(x) \tag{5}
\]
\[
\sum_{N>0}\mathrm{dim}_{\mathbb{Q}}(gr_2^{\mathfrak{D}}\mathcal{H}_N)x^N=\mathbb{O}(x)^2-\mathbb{S}(x)+\mathbb{O}(x)\cdot \mathbb{E}(x)
\]

In the depth $3$ case, we have
\[
gr_3^{\mathfrak{D}}\mathcal{H}\cong gr_3^{\mathfrak{D}}\mathcal{A}\oplus \mathbb{Q}[\zeta^{\mathfrak{m}}(2)]\otimes gr_{2}^{\mathfrak{D}}\mathcal{A}.
\]
So
\[
\sum_{N>0}\mathrm{dim}_{\mathbb{Q}}(gr_3^{\mathfrak{D}}\mathcal{H}_N)x^N=\sum_{N>0}\mathrm{dim}_{\mathbb{Q}}(gr_{3}^{\mathfrak{D}}\mathcal{A}_N)x^N+\sum_{N>0}\mathrm{dim}_{\mathbb{Q}}(gr_{2}^{\mathfrak{D}}\mathcal{A}_N)x^N\cdot\mathbb{E}(x).
\]
By Proposition \ref{half}(i), we have
\[
\begin{split}
\sum_{N>0}\mathrm{dim}_{\mathbb{Q}}(gr_3^{\mathfrak{D}}\mathcal{A}_N)x^N&-\left(\sum_{N>0}\mathrm{dim}_{\mathbb{Q}}(gr_1^{\mathfrak{D}}\mathcal{A}^{odd}_N )x^N \right)\left(\sum_{N>0}\mathrm{dim}_{\mathbb{Q}}(gr_2^{\mathfrak{D}}\mathcal{A}_N)x^N  \right)\\
&+\left(\sum_{N>0}\mathrm{dim}_{\mathbb{Q}}(\mathbb{P}_N^{{}^{\vee}})x^N  \right)\left(\sum_{N>0}\mathrm{dim}_{\mathbb{Q}}(gr_1^{\mathfrak{D}}\mathcal{A})x^N  \right)=0.
\end{split}
\]
So
\[
\sum_{N>0}\mathrm{dim}_{\mathbb{Q}}(gr_3^{\mathfrak{D}}\mathcal{A}_N )x^N)=\mathbb{O}(x)^3-2\mathbb{S}(x)\cdot\mathbb{O}(x) \tag{6}
\]
\[
\sum_{N>0}\mathrm{dim}_{\mathbb{Q}}(gr_3^{\mathfrak{D}}\mathcal{H}_N )x^N=\mathbb{O}(x)^3-2\mathbb{O}(x)\cdot\mathbb{S}(x)+\mathbb{E}(x)\left(\mathbb{O}(x)^2-\mathbb{S}(x)  \right)
\]

For $r\geq 1$, we have
\[
gr_r^{\mathfrak{D}}\mathcal{H}\cong gr_r^{\mathfrak{D}}\mathcal{A}\oplus \mathbb{Q}[\zeta^{\mathfrak{m}}(2)]\otimes gr_{r-1}^{\mathfrak{D}}\mathcal{A}.
\]
So
\[
\sum_{N>0}\mathrm{dim}_{\mathbb{Q}}(gr_r^{\mathfrak{D}}\mathcal{H}_N)x^N=\sum_{N>0}\mathrm{dim}_{\mathbb{Q}}(gr_{r}^{\mathfrak{D}}\mathcal{A}_N)x^N+\sum_{N>0}\mathrm{dim}_{\mathbb{Q}}(gr_{r-1}^{\mathfrak{D}}\mathcal{A}_N)x^N\cdot\mathbb{E}(x).        \tag{7}
\]
For $r\geq 4$, if the long exact sequence conjecture is true, then
\[
\begin{split}
&\left(\sum_{N>0}\mathrm{dim}_{\mathbb{Q}}(\mathbb{P}^{{}^{\vee}}_N)x^N\right)\left(\sum_{N>0}\mathrm{dim}_{\mathbb{Q}}(gr_{r-4}^{\mathfrak{D}}\mathcal{A}_N)x^N  \right)-\sum_{N>0}\mathrm{dim}_{\mathbb{Q}}(gr_r^{\mathfrak{D}}\mathcal{A}_N )x^N\\
&+\left(\sum_{N>0}\mathrm{dim}_{\mathbb{Q}}(gr_1^{\mathfrak{D}}\mathcal{A}^{odd}_N)x^N \right)\left( \sum_{N>0}\mathrm{dim}_{\mathbb{Q}}(gr_{r-1}^{\mathfrak{D}}\mathcal{A}_N)x^N \right)\\
&-\left(\sum_{N>0}\mathrm{dim}_{\mathbb{Q}}(\mathbb{P}^{{}^{\vee}}_N)x^N\right)\left(\sum_{N>0}\mathrm{dim}_{\mathbb{Q}}(gr_{r-2}^{\mathfrak{D}}\mathcal{A}_N)x^N  \right)=0
\end{split}
\]
So
\[
\begin{split}
&\sum_{N>0}\mathrm{dim}_{\mathbb{Q}}(gr_r^{\mathfrak{D}}\mathcal{A}_N )x^N-\mathbb{O}(x)\cdot\sum_{N>0}\mathrm{dim}_{\mathbb{Q}}(gr_{r-1}^{\mathfrak{D}}\mathcal{A}_N )x^N\\
&+\mathbb{S}(x)\cdot\sum_{N>0}\mathrm{dim}_{\mathbb{Q}}(gr_{r-2}^{\mathfrak{D}}\mathcal{A}_N )x^N-\mathbb{S}(x)\cdot\sum_{N>0}\mathrm{dim}_{\mathbb{Q}}(gr_{r-4}^{\mathfrak{D}}\mathcal{A}_N )x^N=0
\end{split}   \tag{8}
\]
From formula $(5)$, $(6)$ and  $(8)$, it follows that
\[
\sum_{N,r>0}\mathrm{dim}_{\mathbb{Q}}(gr_r^{\mathfrak{D}}\mathcal{A}_N  )x^Ny^r=\frac{1}{1-\mathbb{O}(x)y+\mathbb{S}(x)y^2-\mathbb{S}(x)y^4}. \tag{9}
\]
From formula $(7)$ and $(9)$, it follows that
\[
\sum_{N,r>0}\mathrm{dim}_{\mathbb{Q}}(gr_r^{\mathfrak{D}}\mathcal{H}_N  )x^Ny^r=\frac{1+\mathbb{E}(x)}{1-\mathbb{O}(x)y+\mathbb{S}(x)y^2-\mathbb{S}(x)y^4}.
\]
$\hfill\Box$\\

According to the motivic Broadhurst-Kreimer conjecture, we should have \[gr_r^{\mathfrak{D}}\mathcal{A}^{odd}\subsetneqq gr_r^{\mathfrak{D}}\mathcal{A}\]
for $r\geq4$. So in order to prove the conjecture \ref{long}, we need to find the missing generators which should be related to the first map in conjecture \ref{long}.

   From
\[\mathcal{H}\simeq\mathcal{O}(\mathrm{U}^{dR})\otimes\mathbb{Q}[\zeta^{\mathfrak{m}}(2)]\]
we have $\mathcal{A}\cong\mathcal{O}(\mathrm{U}^{dR})$. Since $U^{dR}$ is a pro-unipotent algebraic group, from the long exact sequence conjecture we should have that \[\mathfrak{d}\mathfrak{g}=gr^{\bullet}\mathfrak{g}=\oplus_{r\geq1}\mathfrak{D}^r\mathfrak{g}/\mathfrak{D}^{r+1}\mathfrak{g}\] has generator
$(\mathrm{ad}\,e_0)^{2i}e_1$ in depth $1$ and generator isomorphic to $\mathbb{P}$ in depth $4$ and the only relations between these generators are the period relations between $(\mathrm{ad}\,e_0)^{2i}e_1$. So from the analysis of section $4$ in \cite{enriquez} we know the homology group $H_{n}(\mathfrak{dg},\mathbb{Q}), n\geq 0$ in each weight and depth if the long exact sequence conjecture is true.

We have the following result which may be useful to understand the motivic Broadhurst-Kreimer conjecture.
\begin{Thm}\label{inj}
If for all $r\geq4$, there is an injective map
\[j:\mathbb{P}^{{}^{\vee}}\otimes gr_{r-4}^{\mathfrak{D}}\mathcal{A}\rightarrow gr_r^{\mathfrak{D}}\mathcal{A}
\]
which satisfies $\widetilde{\partial}_g\circ j=0$ and if Conjecture \ref{nondege} is true, then the long exact sequence conjecture is true.
\end{Thm}
\noindent{\bf Proof}:
If such a map $j$ existed for $r\geq4$, assuming Conjecture \ref{nondege}, from Proposition \ref{half}(ii) and formula $(5)$, $(6)$, by the same method in the proof of Theorem \ref{big} we have
\[
(1+\sum_{N,r>0}\mathrm{dim}_{\mathbb{Q}}gr_r^{\mathfrak{D}}\mathcal{A}_Nx^Ny^r)(1-\mathbb{O}(x)y+\mathbb{S}(x)y^2-\mathbb{S}(x)y^4)\geq1.
\]

The inequality means that the coefficient in the series of the left is not smaller than the coefficient  of the right in each corresponding terms.  Furthermore
\[
(1+\sum_{N,r>0}\mathrm{dim}_{\mathbb{Q}}gr_r^{\mathfrak{D}}\mathcal{A}_Nx^Ny^r)(1-\mathbb{O}(x)y+\mathbb{S}(x)y^2-\mathbb{S}(x)y^4)=1
\]
if and only if $\mathrm{Im}\;j=\mathrm{Ker}\;\widetilde{\partial}_g$.

Denote by
\[f(x,y)=(1+\sum_{N,r>0}\mathrm{dim}_{\mathbb{Q}}gr_r^{\mathfrak{D}}\mathcal{A}_Nx^Ny^r)(1-\mathbb{O}(x)y+\mathbb{S}(x)y^2-\mathbb{S}(x)y^4),\]
then we have
\[
\begin{split}
f(x,1)&=(1+\sum_{N,r>0}\mathrm{dim}_{\mathbb{Q}}gr_r^{\mathfrak{D}}\mathcal{A}_Nx^N)(1-\mathbb{O}(x))\\
      &=(1+\sum_{N>0}\mathrm{dim}_{\mathbb{Q}}\mathcal{A}_Nx^N)(1-\mathbb{O}(x))\\
      &=(1+\sum_{N>0}\mathrm{dim}_{\mathbb{Q}}\mathcal{O}(\mathcal{U}^{dR})_N x^N)(1-\mathbb{O}(x))\\
      &=\frac{1}{1-x^3-x^5-\cdots-x^{2n+1}-\cdots}\cdot\left(1-\mathbb{O}(x)\right)\\
      &=1,
\end{split}
\]
where $\mathcal{A}_N$ denote the weight $N$ part of $\mathcal{A}$.

Denote
\[
f(x,y)=1+\sum_{r,s>0}a_{r,s}x^ry^s
\]
Since $f(x,y)\geq1$ and $f(x,1)=1$,  we have
\[a_{r,s}\geq 0,\;\sum_{s>0}a_{r,s}=0,\;\forall r>0.\]
So $\forall r,s>0,\;a_{r,s}=0$.
Then we have $f(x,y)=1$. So the long sequence is exact.$\hfill\Box$\\

In \cite{bro1}, Brown constructed a complex ($(4.1)$ in \cite{bro1})  and conjectured that this complex is an exact sequence (Conjecture $1$ in \cite{bro1}). 
The conjectured exact sequence in \cite{bro1} deals with the formal double shuffle equations and linearized double shuffle Lie algebra $\mathfrak{ls}$ while Conjecture \ref{long} deals with motivic multiple zeta values and depth-graded motivic Lie algebra directly.

\section{Totally odd motivic multiple zeta values in depth two and three}

\subsection{The depth two case}

In this subsection we will consider the polynomial representation of $gr^2_{\mathfrak{D}}\mathcal{U}\mathfrak{g}$. We will show that there are no quasi-uneven elements in $gr^2_{\mathfrak{D}}\mathcal{U}\mathfrak{g}$. As a corollary we obtain a new proof of the fact that $gr_2^{\mathfrak{D}}\mathcal{A}^{odd}=gr_2^{\mathfrak{D}}\mathcal{A}$ which can be proved by Theorem $1$ and Theorem $2$ in \cite{kaneko}.

 By the work of Brown \cite{depth}, the polynomial representation of $gr^2_{\mathfrak{D}}\mathcal{U}\mathfrak{g}$ is the space of  $\mathbb{Q}$-linear combinations of the following polynomials:

 \[
 (y_1-y_0)^{2m_1}(y_2-y_0)^{2m_2}+(y_2-y_1)^{2m_1}(y_1-y_0)^{2m_2}-(y_2-y_1)^{2m_1}(y_2-y_0)^{2m_2}.
 \]
(Hint: Using the formula in subsection \ref{multiplication} to calculate $(y_1-y_0)^{2m_1}\circ(y_1-y_0)^{2m_2},m_1,m_2\geq1$.)

Now we have
\begin{prop}\label{quasiuneven}
For $N\geq6$, denote by $V_{N}$ the $\mathbb{Q}$-vector space which is spanned by the following polynomials:

\[
(y_1-y_0)^{m_1-1}(y_2-y_0)^{m_2-1}+(y_2-y_1)^{m_1-1}(y_1-y_0)^{m_2-1}-(y_2-y_1)^{m_1-1}(y_2-y_0)^{m_2-1},
\]
where $m_1+m_2=N,m_1,m_2\; \mathrm{odd},\geq3.$
Then there are no quasi-uneven elements in $V_N$.
\end{prop}

By definition, Proposition \ref{quasiuneven} is equivalent to the following proposition:
\begin{prop}\label{uneven}
For $N\geq6$, denote by $W_N$ the $\mathbb{Q}$-vector space which is spanned by the following elements

\[
y_1^{m_1-1}y_2^{m_2-1}+(y_2-y_1)^{m_1-1}(y_1^{m_2-1}-y_2^{m_2-1}),m_1+m_2=N,m_1,m_2\geq3,\mathrm{odd}.
\]
If $f\in W_N$ and the coefficients of $y_1^{n_1-1}y_2^{n_2-1},n_1,n_2\geq3,\;odd$ in $f$ are all zero, then $f=0$.
\end{prop}
\noindent{\bf Proof}:
Denote by  \[f=\sum_{\overrightarrow{m}\in S_{N,2}}p_{m_1,m_2}[y_1^{m_1-1}y_2^{m_2-1}+(y_2-y_1)^{m_1-1}(y_1^{m_2-1}-y_2^{m_2-1})]\]
and
\[
p=\sum_{\overrightarrow{m}\in S_{N,2}}p_{m_1,m_2}y_1^{m_1-1}y_2^{m_2-1},
\]
then \[f(y_1,y_2)=p(y_1,y_2)+p(y_2-y_1,y_1)-p(y_2-y_1,y_2).\]

By explicit calculation, the totally even part of $$y_1^{m_1-1}y_2^{m_2-1}+(y_2-y_1)^{m_1-1}(y_1^{m_2-1}-y_2^{m_2-1})$$
is
\[
\sum_{\overrightarrow{n}\in S_{N,2}}e\dbinom{m_1,m_2}{n_1,n_2}y_1^{n_1-1}y_2^{n_2-1}
\]

If the coefficients of $y_1^{n_1-1}y_2^{n_2-1},n_1,n_2\geq3,\mathrm{odd}$ in $f$ are all zero,
 then we have $\pi_1(p)\in \mathrm{Ker}\;E_{N,2}$. Since  $\pi_1(W_{N,2})=\mathrm{Ker}\;E_{N,2}$ (main results in Baumard, Schneps \cite{schneps} which is reproved in Tasaka\cite{tasaka}), we deduce that $p$ is an even restricted period polynomial, i.e. $f=0$.
 $\hfill\Box$\\

From Proposition \ref{quasiuneven} we know that there are no quasi-uneven elements in $gr_{\mathfrak{D}}^2\mathcal{U}\mathfrak{g}$. As a corollary we have:

\begin{Cor}\label{newproof}
$gr_{\mathfrak{D}}^2\mathcal{A}^{odd}=gr_{\mathfrak{D}}^2\mathcal{A}$.
\end{Cor}

\subsection{The depth three case}
By Theorem $2.6$ in Goncharov \cite{goncharov}, the natural map $\beta:\mathrm{Lie}_n(\mathfrak{d}\mathfrak{g}_1)\rightarrow \mathfrak{d}\mathfrak{g}_n$ is surjective for $n=3$. So in $\mathcal{U}\mathfrak{h}$, $gr_{\mathfrak{D}}^3\mathcal{U}\mathfrak{g}$ is generated by elements of the form $$(\mathrm{ad}\,e_0)^{2n_1}e_1\otimes(\mathrm{ad}\,e_0)^{2n_2}e_1\otimes(\mathrm{ad}\,e_0)^{2n_3}e_1.$$

By the definition of $\circ$ we have
\[
\begin{split}
&\;\;\;\;\gamma((\mathrm{ad}\,e_0)^{2n_1}e_1\otimes(\mathrm{ad}\,e_0)^{2n_2}e_1\otimes(\mathrm{ad}\,e_0)^{2n_3}e_1)\\
&=(\mathrm{ad}\,e_0)^{2n_1}e_1\circ(\mathrm{ad}\,e_0)^{2n_2}e_1\circ(\mathrm{ad}\,e_0)^{2n_3}e_1\\
&=(\mathrm{ad}\,e_0)^{2n_1}e_1\circ\left((\mathrm{ad}\,e_0)^{2n_2}e_1\circ(\mathrm{ad}\,e_0)^{2n_3}e_1\right)
\end{split}
\]
By Proposition $2.2$, Subsection $6.1$ in \cite{depth} and Proposition \ref{action} in Section \ref{key}, the polynomial representation of $(\mathrm{ad}\,e_0)^{2n_1}e_1\circ\left((\mathrm{ad}\,e_0)^{2n_2}e_1\circ(\mathrm{ad}\,e_0)^{2n_3}e_1\right)$
is
\[
\begin{split}
&(y_1-y_0)^{2m_1}[(y_2-y_0)^{2m_2}(y_3-y_0)^{2m_3}+(y_3-y_2)^{2m_2}(y_2-y_0)^{2m_3}\\
&\;\;\;\;\;\;\;\;\;\;\;\;\;\;\;\;\;\;\;\;\;\;\;\;\;\;\;\;\;\;\;\;\;\;\;\;\;\;\;\;\;\;\;\;\;\;\;\;\;\;\;\;\;\;\;-(y_3-y_2)^{2m_2}(y_3-y_0)^{2m_3}]
\end{split}
\]
\[
\begin{split}
&+(y_2-y_1)^{2m_1}[(y_1-y_0)^{2m_2}(y_3-y_0)^{2m_3}+(y_3-y_1)^{2m_2}(y_1-y_0)^{2m_3}\\
&\;\;\;\;\;\;\;\;\;\;\;\;\;\;\;\;\;\;\;\;\;\;\;\;\;\;\;\;\;\;\;\;\;\;\;\;\;\;\;\;\;\;\;\;\;\;\;\;\;\;\;\;\;\;\;\;\;-(y_3-y_1)^{2m_2}(y_3-y_0)^{2m_3}]
\end{split}
\]
\[
\begin{split}
&+(y_3-y_2)^{2m_1}[(y_1-y_0)^{2m_2}(y_2-y_0)^{2m_3}+(y_2-y_1)^{2m_2}(y_1-y_0)^{2m_3}\\
&\;\;\;\;\;\;\;\;\;\;\;\;\;\;\;\;\;\;\;\;\;\;\;\;\;\;\;\;\;\;\;\;\;\;\;\;\;\;\;\;\;\;\;\;\;\;\;\;\;\;\;\;\;\;\;\;\;-(y_2-y_1)^{2m_2}(y_2-y_0)^{2m_3}]
\end{split}
\]
\[
\begin{split}
&-(y_2-y_1)^{2m_1}[(y_2-y_0)^{2m_2}(y_3-y_0)^{2m_3}+(y_3-y_2)^{2m_2}(y_2-y_0)^{2m_3}\\
&\;\;\;\;\;\;\;\;\;\;\;\;\;\;\;\;\;\;\;\;\;\;\;\;\;\;\;\;\;\;\;\;\;\;\;\;\;\;\;\;\;\;\;\;\;\;\;\;\;\;\;\;\;\;\;\;\;-(y_3-y_2)^{2m_2}(y_3-y_0)^{2m_3}]
\end{split}
\]
\[
\begin{split}
&-(y_3-y_2)^{2m_1}[(y_1-y_0)^{2m_2}(y_3-y_0)^{2m_3}+(y_3-y_1)^{2m_2}(y_1-y_0)^{2m_3}\\
&\;\;\;\;\;\;\;\;\;\;\;\;\;\;\;\;\;\;\;\;\;\;\;\;\;\;\;\;\;\;\;\;\;\;\;\;\;\;\;\;\;\;\;\;\;\;\;\;\;\;\;\;\;\;\;\;\;-(y_3-y_1)^{2m_2}(y_3-y_0)^{2m_3}]
\end{split}
\]

In the above polynomial, the part which correspond to $y_1^{m_1-1}y_2^{m_2-1}y_3^{m_3-1},n_i,\geq 3,\mathrm{odd},1\leq i\leq 3$ is
\[
\sum_{\overrightarrow{n}\in S_{N,3}}c\dbinom{m_1,m_2,m_3}{n_1,n_2,n_3}y_1^{n_1-1}y_2^{n_2-1}y_3^{n_3-1}.
\]
Where \[
C_{N,3}=\left(c\dbinom{m_1,m_2,m_3}{n_1,n_2,n_3}\right)
_{\substack {\overrightarrow{m}\in S_{N,3}\\\overrightarrow{n}\in S_{N,3}}}=E^{(2)}_{N,3}\cdot E_{N,3}.\]

Denote by $(gr^3_{\mathfrak{D}}\mathcal{U}\mathfrak{g})_N$ the weight $N$ part of $gr^3_{\mathfrak{D}}\mathcal{U}\mathfrak{g}$, by Proposition \ref{half}(i), we have
\[
\sum_{N>0}\mathrm{dim}_{\mathbb{Q}}(gr^3_{\mathfrak{D}}\mathcal{U}\mathfrak{g})_N x^N=\mathbb{O}(x)^3-2\mathbb{O}(x)\mathbb{S}(x).
\]

So we have $gr^{\mathfrak{D}}_3\mathcal{A}^{odd}=gr^{\mathfrak{D}}_3\mathcal{A}\Leftrightarrow$ there are no quasi-uneven elements in $gr^3_{\mathfrak{D}}\mathcal{U}\mathfrak{g}\Leftrightarrow$
$\sum_{N>0}\mathrm{rank}\; C_{N,3}x^N=\mathbb{O}(x)^3-2\mathbb{O}(x)\mathbb{S}(x)$.

By the proof of Proposition \ref{sur} and the results of Tasaka, we know the map
\[\xi:\mathbf{W}_{N,3}\longrightarrow \mathbf{Vect}_{N,3}\]
\[p(x_1,x_2,x_3)\longmapsto (\pi_1(p(x_1,x_2,x_3))E^{(2)}_{N,3}\]
satisfies $\mathrm{Im}\;\xi\subseteq\mathrm{Ker} \;E_{N,3}$. By Goncharov's theorem $\mathbb{P}\otimes\mathfrak{d}\mathfrak{g}_1\cap\wedge^3\mathfrak{d}\mathfrak{g}_1={0}$, $\xi$ is injective. So $\xi$ induces an injective map
\[\widetilde{\xi}:\mathbf{W}_{N,3}\longrightarrow \mathrm{Ker}\; E_{N,3}\]
\[p(x_1,x_2,x_3)\longmapsto (\pi_1(p(x_1,x_2,x_3))E^{(2)}_{N,3}.\]

The following conjecture is due to Tasaka.
\begin{conj}\label{iso}
$\widetilde{\xi}$ is isomorphic.
\end{conj}
Now we have
\begin{prop}\label{three}
If Conjecture \ref{iso} is true, then $gr_3^{\mathfrak{D}}\mathcal{A}^{odd}=gr_3^{\mathfrak{D}}\mathcal{A}$.
\end{prop}
\noindent{\bf Proof}:
Since $C_{N,3}=E^{(2)}_{N,3}\cdot E_{N,3}$, so
\[
\mathrm{Ker}\;C_{N,3}=\mathrm{Ker}\;E_{N,3}^{(2)}\oplus\mathrm{Im}\;E_{N,3}^{(2)}\cap\mathrm{Ker}\;E_{N,3}.
\]

If $\widetilde{\xi}$ is isomorphic, we will have $\mathrm{Ker}\;E_{N,3}\subseteq\mathrm{Im}\;E_{N,3}^{(2)}$ and $\mathrm{Ker}\; E_{N,3}\cong W_{N,3}$.
So
 \[
\mathrm{Ker}\;C_{N,3}=\mathrm{Ker}\;E_{N,3}^{(2)}\oplus W_{N,3}
\]
and \[\sum_{N>0}\mathrm{rank}\; C_{N,3}x^N=\mathbb{O}(x)^3-2\mathbb{O}(x)\mathbb{S}(x).\]
From the above formula we have $gr_3^{\mathfrak{D}}\mathcal{A}^{odd}=gr_3^{\mathfrak{D}}\mathcal{A}$.$\hfill\Box$\\
\begin{rem}
The proof of the injectivity of $\xi$  is indirect and based on the result $\mathbb{P}\otimes\mathfrak{d}\mathfrak{g}_1\cap\wedge^3\mathfrak{d}\mathfrak{g}_1={0}$. A direct proof of the injectivity of $\xi$ should be useful to tackle Conjecture \ref{iso}.
\end{rem}

\section*{Acknowledgement}
In this paper, the author wants to thank Pierre Deligne for his explanations of the author's questions about mixed Tate motives. The author also thanks Yuancao Zhang for his help about understanding many details during the author's study. The author also thanks Leila Schneps for her detailed explanation of a fact in her and Samuel Baumard's paper \cite{schneps}. At last, the author wants to thank his supervisor Qingchun Tian for his suggestion to choose this subject and helpful advice during the research.


\newpage

\begin{thebibliography}{1}


\bibitem{schneps}
   S. Baumard, L. Schneps,
    \emph{Period polynomial relations between double zeta values},
    Ramanujian J.,
    32(1) (2013),
    83-100.

\bibitem{kreimer}
   D. Broadhurst, D. Kreimer,
   \emph{Association of multiple zeta values with positive knots via Feynman diagrams up to $9$ loops},
   Phys. Lett. B 393 (1997), no. 3-4, 303-412.

\bibitem{bro1}
   F. Brown,
   \emph{An exact sequence for the Broadhurst-Kreimer conjecture},
   http://www. ihes. fr/~brown /BKExactSeq1. pdf.

\bibitem{depth}
   F. Brown,
   \emph{Depth-graded motivic multiple zeta value},
   arXiv: 1301. 3053.

\bibitem{brown}
     F. Brown,
     \emph{Mixed Tate motives over $\mathbb{Z}$},
     Ann. of Math.,
     175(2) (2012),
     949-976.

\bibitem{zeta3}
    F. Brown,
    \emph{Zeta elements in depth $3$ and the fundamental Lie algebra of a punctured elliptic curve},
    arXiv: 1504. 04737.

\bibitem{BF}
  J. I. Burgos Gil,  J. Fresán,
 \emph{Multiple zeta values: from numbers to motives},
 Clay Math. Proceedings, to appear
\bibitem{cartier}
   P. Cartier,
   \emph{A primer of hopf algebras},
   Frontiers in Number Theory,
   Physics, and Geometry II, (2007),
   537-615.

\bibitem{deligne}
   P. Deligne, A. B. Goncharov,
   \emph{Groupes fondamentaux motiviques de Tate mixte},
    Ann. Sci. \'{E}cole Norm. Sup.,
    38 (2005),
    1-56.

\bibitem{enriquez}
   B. Enriquez, P. Lochak,
   \emph{Homology of depth-graded motivic Lie algebras and koszulity},
   arXiv: 1407. 4060.


\bibitem{kaneko}
    H. Gangl, M. Kaneko, D. Zagier,
    \emph{Double zeta values and modular forms},
    Automorphic forms and zeta functions,
    In: Proceedings of the conference in memory of Tsuneo Arakawa,
    World Scientific (2006),
    71-106.

\bibitem{goncha}
    A. B. Goncharov,
    \emph{Galois symmetries of fundamental groupoids and noncommutative goemetry},
    Duke Math. J., 128 (2) (2005), 209-284.

\bibitem{goncharov}
    A. B. Goncharov,
    \emph{Multiple polylogarithms, cyclotomy and modular complexes},
    Mathematical Research Letters 5 (1998), 497-516.
\bibitem{gon}
    A. B. Goncharov, \emph{The dihedral Lie algebras and Galois symmetries of $\pi_1^{(l)}(\mathbb{P}^1-\{0,\infty\}\cup \mu_{N})$},
    Duke Math. J., 110 (3)  (2001), 397-487.

\bibitem{hoff}
    M. E. Hoffman,
    \emph{The algebra of multiple harmonic series},
    J.Algebra, 194(2) (1997), 477-495.

\bibitem{ihara}
    K. Ihara, M. Kaneko, D, Zagier,
    \emph{Derivation and double shuffle relations for multiple zeta values},
    Compos. Math, 142 (2006), 307-308.

\bibitem{ihara1}
   Y. Ihara,
   \emph{Some arithmetic aspects of Galois actions on the pro-p fundamental group of  $\mathbb{P}-\{0,1,\infty\}$},
   in:  Arithmetic Fundamental Groups
and Noncommutative Algebra, Proc. Sympos. Pure Math. 70,  Berkeley,
CA, (1999,) 247-273.

\bibitem{ll}
   J. Li, F. Liu,
   \emph{Motivic double zeta values of odd weight},
   arXiv: 1710. 02244.
\bibitem{racinet}

    G. Racinet, \emph{Doubles m\'{e}langes des polylogarithmes multiples aux racines de l’unit\'{e}},
    Publ. Math. Inst. Hautes \'{E}tudes Sci.  95 (2002), 185-231.


\bibitem{sch}
    L. Schneps,
    \emph{On the Poisson Bracket on the Free Lie Algebra in two Generators},
    Journal of Lie Theory
    16, no. 1 (2006),
    19-37.
\bibitem{souderes}
    I. Soud\`{e}res,
    \emph{Motivic double shuffle},
    Int. J . Number Theory 6 (2010),
    339-370.
\bibitem{tasaka}
   K. Tasaka,
   \emph{On linear relations among totally odd multiple zeta values related to period polynomials},
   arXiv: 1402. 3391.








\end{thebibliography}
\end{document}